%% file: arxivv1.tex
\documentclass[onefignum,onetabnum]{siamart250211}

\input{preamble}
\input{colordef}
\input{macros}
\input{siam-shared}

\ifpdf
\hypersetup{
  pdftitle={A Performance Portable Fast Ewald Summation for Stokes Flow},
  pdfauthor={G. Kosmacher, J. Bagge, Z. Du, G. Biros}
}
\fi

\begin{document}

\maketitle

\begin{abstract}
	\input{sections/abstract}
\end{abstract}



% Paper Body
\input{sections/introduction}

\input{sections/methods}

\input{sections/methods_gpu}

\input{sections/methods_erf}
\input{sections/methods_model}
\input{sections/results}

\input{sections/conclusion}

\section*{Acknowledgments}
The authors acknowledge the Texas Advanced Computing Center (TACC)\footnote{URL: \url{http://www.tacc.utexas.edu}} at The University of Texas at Austin and the Oregon Advanced Computing Institute for Science and Society (OACISS)\footnote{URL: \url{https://blogs.uoregon.edu/oaciss/}} at The University of Oregon for providing HPC resources that have contributed to the research results reported within this paper. JB was supported in part by the Peter O'Donnell Jr.\ Postdoctoral Fellowship at the Oden Institute.

\bibliographystyle{siamplain}
\bibliography{ref.bib}

\input{sections/appendix}

\end{document}

%% file: preamble.tex
\usepackage{amsmath}
\usepackage[inline]{enumitem}
\usepackage{nicefrac}
\usepackage{enotez}
\usepackage{booktabs, siunitx, colortbl, multirow}
\usepackage{tikz}
\usepackage{pythonhighlight}
\usepackage{mathtools}

\usepackage{algpseudocode}

\algrenewcommand\algorithmiccomment[1]{\hfill\textcolor{gray}{// #1}}

\algblock{ParFor}{EndParFor}
\algnewcommand\algorithmicparfor{\textbf{parfor}}
\algnewcommand\algorithmicpardo{\textbf{do}}
\algnewcommand\algorithmicendparfor{\textbf{end\ parfor}}
\algrenewtext{ParFor}[1]{\algorithmicparfor\ #1\ \algorithmicpardo}
\algrenewtext{EndParFor}{\algorithmicendparfor}

%% file: colordef.tex
\usepackage{xcolor}

\definecolor{orange}{HTML}{f8971f}
\definecolor{yellow}{HTML}{ffd600}
\definecolor{lightgreen}{HTML}{a6cd57}
\definecolor{green}{HTML}{579d42}
\definecolor{lightblue}{HTML}{00a9b7}
\definecolor{blue}{HTML}{005f86}
\definecolor{grey}{HTML}{9cadb7}
\definecolor{offwhite}{HTML}{d6d2c4}
\definecolor{dark}{HTML}{333f48}

%% file: macros.tex
\newcommand{\n}{\boldsymbol{n}}
\newcommand{\g}{\boldsymbol{g}}
\newcommand{\p}{\boldsymbol{p}}
\newcommand{\q}{\boldsymbol{q}}
\newcommand{\x}{\boldsymbol{x}}
\newcommand{\y}{\boldsymbol{y}}
\newcommand{\f}{\boldsymbol{f}}
\newcommand{\G}{\boldsymbol{G}}
\newcommand{\T}{\boldsymbol{T}}

\renewcommand{\u}{\boldsymbol{u}}
\renewcommand{\r}{\boldsymbol{r}}
\renewcommand{\a}{\boldsymbol{a}}
\newcommand{\zv}{\boldsymbol{0}}
\newcommand{\B}{\mathcal{B}}
\newcommand{\E}{\mathrm{e}}

\newcommand{\I}{\mathrm{i}}
\newcommand{\Id}{\boldsymbol{I}}
\newcommand{\erf}{\mathrm{erf}}

\newcommand{\K}{\mathcal{K}}
\newcommand{\rc}{r_{\mathrm{c}}}
\newcommand{\Ns}{N_{\mathrm{s}}}
\newcommand{\Nt}{N_{\mathrm{t}}}
\newcommand{\Ng}{N_{\mathrm{g}}}
\newcommand{\NGPU}{N_{\mathrm{GPU}}}
\newcommand{\Ps}{\mathcal{P}}
\newcommand{\N}{\mathrm{N}}
\newcommand{\F}{\mathrm{F}}
\newcommand{\Cell}{\ensuremath{\mathcal{C}}}
\newcommand{\CF}{\ensuremath{\mathcal{C}^\F}}
\newcommand{\CN}{\ensuremath{\mathcal{C}^\N}}
\newcommand{\GN}{\G^{\N}}
\newcommand{\GF}{\G^{\F}}
\DeclareMathOperator{\erfc}{erfc}
\newcommand{\Fou}{\mathcal{F}_1}
\newcommand{\Gh}{\widehat{\G}}
\newcommand{\kv}{\boldsymbol{\kappa}}
\newcommand{\eSE}{\varepsilon}
\newcommand{\pylib}{\texttt{ParkiPy}}

\newcommand{\U}{\textbf{U}\,}
\renewcommand{\S}{\textbf{S}\,}
\renewcommand{\N}{\textbf{N}\,}

%% file: siam-shared.tex
\usepackage{listings}
\usepackage{lipsum}
\usepackage{amsfonts}
\usepackage{graphicx}
\usepackage{epstopdf}
\ifpdf
  \DeclareGraphicsExtensions{.eps,.pdf,.png,.jpg}
\else
  \DeclareGraphicsExtensions{.eps}
\fi

\newsiamremark{remark}{Remark}
\newsiamremark{hypothesis}{Hypothesis}
\crefname{hypothesis}{Hypothesis}{Hypotheses}
\newsiamthm{claim}{Claim}
\newsiamremark{fact}{Fact}
\crefname{fact}{Fact}{Facts}

\headers{Fast Ewald Summation for Stokes Flow}{G. Kosmacher, Z. Du, and J. Bagge, G. Biros}

\title{A Performance Portable Fast Ewald Summation for Stokes Flow}

\author{Gabriel Kosmacher\thanks{Oden Institute for Computational Engineering and Sciences, The University of Texas, Austin, TX 
  (\email{\{gkosmacher, ziyu\_du\}@utexas.edu}, joar.bagge@austin.utexas.edu, biros@oden.utexas.edu).}
\and Ziyu Du\footnotemark[2]
\and Joar Bagge\footnotemark[2]
\and George Biros\footnotemark[2]}

\usepackage{amsopn}

%% file: sections/abstract.tex
We present GPU algorithms for Ewald summation methods for accelerating N-body Stokes flow problems in periodic domains. Like most N-body codes, Ewald sums use a near-field/far-field decomposition. The near field involves particle-to-particle (P2P) interactions. 
The far field primarily involves particle-to-grid (P2G) and grid-to-particle (G2P) interactions, as well as Fast Fourier Transforms.  
For each interaction, we investigate several algorithmic variants.
Our implementation uses PyKokkos, a Python interface for  the Kokkos C++ parallel programming framework, which supports portability to AMD/\mbox{NVIDIA} GPU and ARM/x86 CPU architectures. 
Double and single-precision numerical results, alongside analytical performance models, confirm the efficiency of our algorithms on AMD and NVIDIA GPU and on ARM and AMD CPU architectures. The P2P interaction achieves around 73\% compute efficiency on NVIDIA H200, 84\% on NVIDIA A100,  60\% on AMD MI300, 52\% on Grace CPU, and 68\% on AMD Epyc CPU.
A straightforward implementation of the P2G kernel can become a computational bottleneck. We introduce a novel P2G algorithm that achieves up to 16$\times$ speedup compared to a baseline GPU implementation.
The overall Ewald sum code processes approximately 8 million particles per second on a H200 GPU, and about a half-million particles per second on a Grace CPU, for nine digits of accuracy. We also perform a multi-GPU weak scaling test on up to 256 million particles (64 GPUs) that shows bounded communication cost for all stages except the all-to-all particle sorting, which can be reduced to neighbor communication in the relevant time-stepping regime. 

%% file: sections/introduction.tex
\section{Introduction}

Stokes suspensions involve the interaction of viscous flows with rigid or deformable particles~\cite{yan-malhotra-shelley20}.  They appear in microfluidic devices~\cite{holm-fedosov-e19}, swimming microorganisms~\cite{Guasto_2012,Ladiges_2021}, capillary blood flow~\cite{Squires_2005}, antibodies~\cite{Lai_2021},  and industrial emulsions~\cite{zinchenko-davis03}. 
Boundary integral equations are a commonly used numerical method for simulating Stokes suspensions~\cite{pozrikidis-00}.
Upon discretization, the velocity field (henceforth referred to as \emph{potential}  per N-body literature conventions)  is written as a sum over $\Ns$ \emph{source  points} $\y_j$, of the form 
\begin{equation}
\label{eq:N-body-sum}
\u(\x_i) = \sum_{j=1}^{\Ns} \G(\x_i-\y_j) \f(\y_j), \quad i =1,\ldots, \Nt,
\end{equation}
where $\G$ is a Stokes  Green's function and $\f(\y_j)$ is a known force density (or simply \emph{density}) associated with $\y_j$.
The sum in \cref{eq:N-body-sum} needs to be evaluated at $\Nt$ \emph{target points} $\x_i$;
assuming $N=\Nt=\Ns$,  direct evaluation costs $\mathcal{O}(N^2)$.
N-body methods reduce this cost to $\mathcal{O}(N)$ for the fast multipole method (FMM) \cite{wang2021exafmm,fmm,lashuk-biros-e09} and to  $\mathcal{O}(N\log N)$ for an Ewald  method~\cite{Lindbo_2010}. 
Although the FMM has favorable asymptotic complexity, the prefactor constants can be quite high compared to the Ewald method when the spatial particle distribution is nearly uniform~\cite{gholami-malhotra-sundar-biros16}.
As a result, Ewald sum methods are widely used for dense suspensions with periodic boundary conditions~\cite{Bagge_2021, Wang_2016, Klinteberg_2016, zhao-freund10}.
\emph{Yet, as far as we are aware, there is no existing open source GPU implementation for a Stokes Ewald method.}

\pagebreak
{\bf Contributions:} Following the Stokes Ewald algorithm introduced in \cite{Bagge_Tornberg_2023} (see \cref{sec:sequential}), 
we present the following methodological and experimental results:
\begin{itemize}[topsep=0pt,itemsep=-1ex,partopsep=1ex,parsep=1ex]

\item The design and implementation of GPU algorithms for different steps in Ewald sums; in particular, our novel particle-to-grid (P2G) method achieves a speedup of more than 16$\times$ compared to a baseline implementation (see \cref{sec:acceleration} and \cref{sec:gpu-p2g}).
\item A novel rational function approximation for $\erf(x)/x$ that appears in Ewald sums; it is faster than system libraries in lower precision (see \cref{sec:special}).
\item A roofline analysis for our kernels with various assumptions that aligns with the measured times, enabling automatic selection of algorithm parameters (see \cref{sec:perf-mod}).
\item Evaluation on multiple platforms, including NVIDIA H200 and A100 GPUs, AMD MI300A GPU, AMD Epyc CPU,  and NVIDIA Grace ARM CPU (see \cref{sec:results}).
\item An open source Python library \pylib{}\footnote{\url{https://github.com/ut-padas/parki}} (supplement \ref{sec:parkipy}), developed using PyKokkos, supporting single node APIs for Stokes and Poisson potentials in single and double precision and a distributed API for the Stokes potential in double precision. 
\end{itemize}
Using standard techniques, we also present results using the message-passing interface (MPI) on up to 64 GPUs.

{\bf Related work.}
Ewald methods, also known as smooth particle-mesh Ewald sums,  have been extensively studied for the Poisson problem for molecular dynamics simulations.
The Ewald decomposition for Stokes was introduced by Hasimoto~\cite{hasimoto1959periodic}.
As mentioned, our summation algorithm is based on the fast spectrally accurate Ewald summation method in \cite{Bagge_Tornberg_2023}, which is a refinement of the scheme introduced in \cite{Lindbo_2010}.  Spectral Ewald summations are similar to smoothed particle mesh methods but more accurate~\cite{Lindbo_2010,zhao-freund10,shamshirgar2021fast}. 

The development of parallel Ewald sums has been largely focused on Poisson kernels for molecular dynamics~\cite{harvey2009spme,phillips2014namd,anton1,anton3,sbalzarini2006ppm}. Examples of work on GPUs include GROMACS~\cite{gromacs}, NAMD~\cite{phillips2020scalable},  AMBER~\cite{salomon2013routine} and many others~\cite{jasz2020classical}.  The primary objective in most studies is to assess the overall performance of an atomistic simulation in terms of nanoseconds per day. Surprisingly, detailed discussions and performance evaluations for the Ewald sum and its components are rare compared to the richer FMM literature.
For Stokes flows, implementations of GPU-accelerated Ewald methods are mentioned in Wang and Brady \cite{Wang_2016} and Fiore et al. \cite{Fiore_2017}. These are restricted to smaller systems for Brownian dynamics (fewer than 10K particles), with no assessment of the Ewald summation component. The work on GPU nonuniform fast Fourier transforms~\cite{cuFINUFFT} is one of the few exceptions that examines different compute steps that closely resemble parts of Ewald sum. 

In terms of portability, frameworks like OpenCL \cite{opencl}, Kokkos \cite{Kokkos0,Kokkos3}, and RAJA \cite{raja} have been developed to support code execution on CPUs and GPUs. Kokkos has been used to develop Poisson Ewald codes in LAMMPS \cite{lammps, Halver}. Another portable implementation is from Mayani et al.\ \cite{Kokkos_Poisson}, who implemented a free-space spectral Poisson solver via Kokkos. PyKokkos \cite{pyKokkos,KokkosEco} enables programming Kokkos in Python with minimal overheads using just-in-time compiling.

%% file: sections/methods.tex
\section{Sequential Ewald Summation}
\label{sec:sequential}

\begin{figure}
  \centering
  \raisebox{15pt}{\includegraphics[width=0.6\textwidth, trim=5mm 5mm 5mm 5mm, clip]{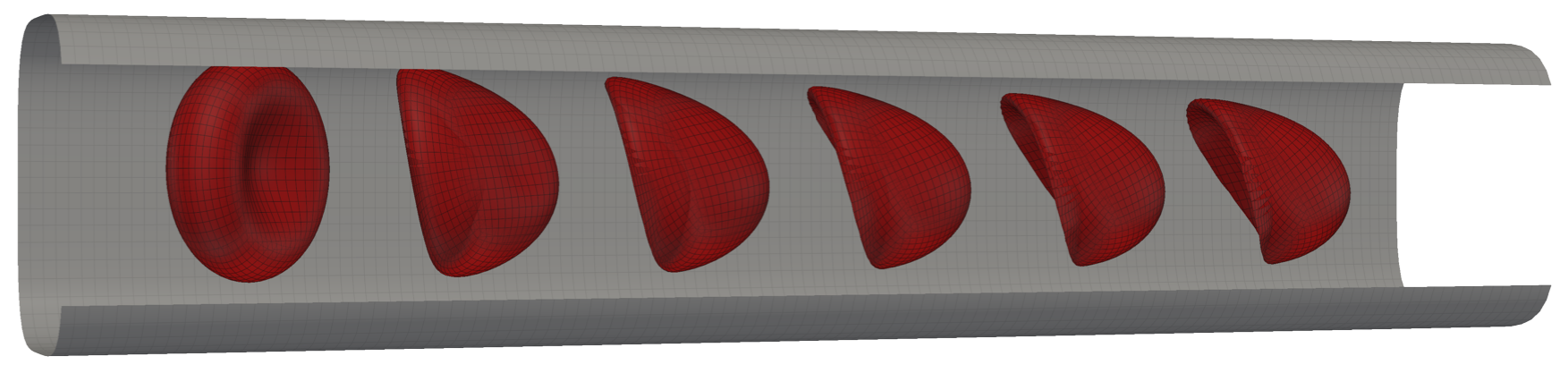}}%
  \hspace{25pt}%
  \includegraphics[width=0.3\textwidth]{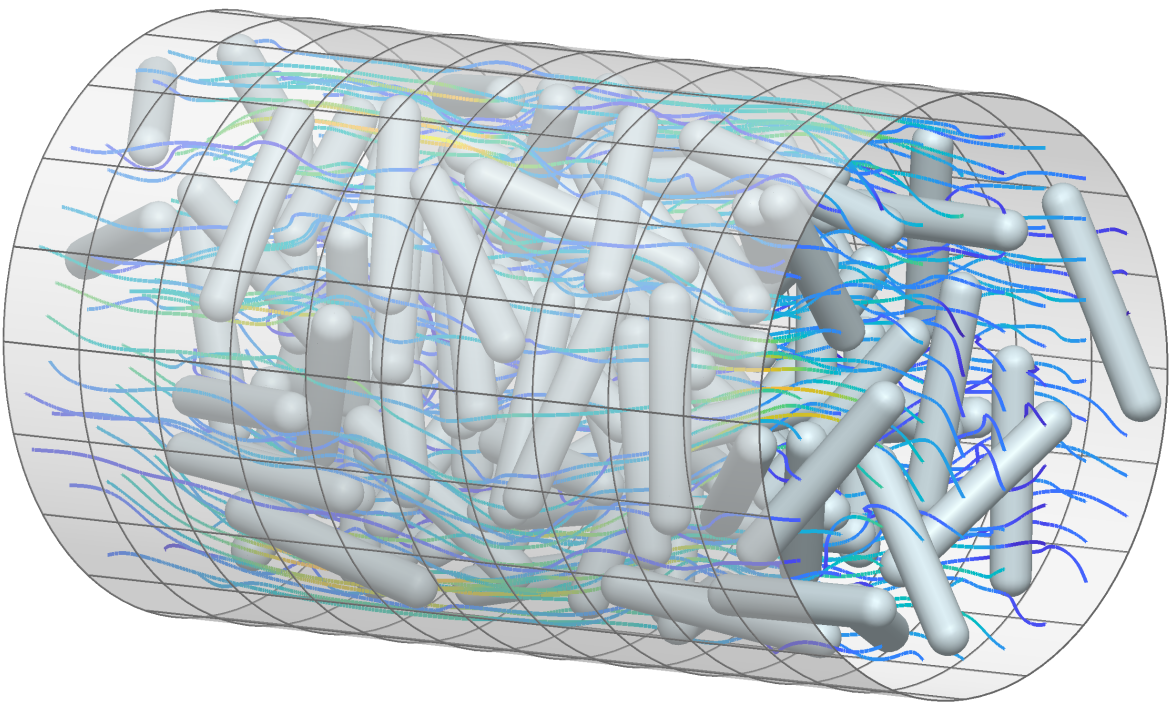}%
  \caption{Two examples of Stokesian flows. \emph{Left}:  A simulation from ongoing work in our group showing deformable particles (red blood cells) flowing in a micrometer-wide pipe. \emph{Right}:  A suspension of a large number of rigid particles (slender rods) in a pipe, taken from~\cite{Bagge_2021}. 
  Both problems are periodic in the direction parallel to the axis of the pipe. Each time step involves particle-related calculations followed by a Stokes Ewald sum.}
  \label{fig:teaser}
\end{figure}

We summarize the Ewald sum scheme from~\cite{Bagge_Tornberg_2023}.
The periodized  version of \cref{eq:N-body-sum} reads
\begin{equation}
\label{eq:potentials}
    \u(\x_i) = \sum_{j=1}^{\Ns} \sum_{\p \in \Ps} \G(\x_i - \y_j + \p) \f(\y_j),
    \qquad i=1,\ldots,\Nt,
\end{equation}
where $\Ps$ denotes a set of periodic images that tile $\mathbb{R}^d$ to create a periodic field. Depending on the problem, we can have periodicity in one, two, or all three directions.
Here for simplicity we only consider the 1-periodic case, $\Ps = \{ (L_1 \alpha, 0, 0) : \alpha \in \mathbb{Z} \}$, with the other directions being free;
Ewald sums of this case can be used to accelerate the solution of confined Stokes solvers in pipe geometries, such as those in \cref{fig:teaser}.
The extension to the 2- and 3-periodic cases is straightforward and does not require any changes to the kernels discussed in \cref{sec:acceleration}.
We assume $\{\x_i\}$, $\{\y_j\}$ $\subset \B \coloneqq [0,L_1) \times [0,L_2) \times [0,L_3)$, and that $\B$ is periodically replicated along the $x_1$ axis.
The Green's function  $\G$ is the Stokeslet given by
$
\G(\r) = \frac{\Id}{\lVert \r \rVert} + \frac{\r \otimes \r}{\lVert \r \rVert^3},
$
where $\Id$ is the 3-by-3 identity matrix. Note that if $\|\r\|\approx \zv$, $\G$ is singular; if $\y_j\approx \x_i$, the term is skipped.
Ewald summation splits $\u = \u^\N + \u^\F$ by setting $\G = \GN + \GF$. 
The near-field component $\GN$ rapidly decays with $\|\r\|$; the far-field component $\GF$ is smooth and slowly decays with $\|\r\|$. 
A \emph{split parameter} $\xi>0$ controls the decay of $\GN$.

\subsection{\texorpdfstring{$\GN$}{GN}: Near-Field (P2P)} \label{sec:sequential-near}
The near-field sum (or P2P interaction) is given by
\begin{equation}
\label{eq:potentials-near}
    \u^\N(\x_i) = \sum_{j=1}^{\Ns} \sum_{\p \in \Ps} \GN(\x_i - \y_j + \p) \cdot \f(\y_j),
\end{equation}
and the expression for the near-field Stokeslet $\GN$ is
\begin{equation}\label{eq:g-near}
  \GN(\r) =
 \G(\r) 
\left(
\erfc(\xi \lVert \r \rVert)
+ \frac{\|\r\| 2\xi \E^{-\xi^2\lVert \r \rVert^2}}{\sqrt{\pi}}
\right)
- \Id \frac{4\xi \E^{-\xi^2\lVert \r \rVert^2}}{\sqrt{\pi}},
\end{equation}
where $\erfc(x) = 1-\erf(x)$ is the complementary error function.
Let the \emph{cutoff radius} $\rc \sim \xi^{-1}$ be a function of the split parameter $\xi$ whose form is derived via \cref{eq:ese}.
$\G^\N$ can then be efficiently computed by dropping any interaction with $\|\r\|>\rc$.
To exploit this locality in \cref{eq:g-near}, we spatially partition $\B$ into uniform cubic \emph{cells} of width $\rc$ and then assign each source and target point to a cell.
When evaluating \cref{eq:potentials-near} on a target $\x_i$, it is enough to determine which cell $\x_i$ belongs to and then restrict its P2P interactions to the 27-cell neighborhood of that cell.
To control the near-field error, we need to choose $\xi$ appropriately (see \cref{sec:sequential-accuracy}).  

\subsection{\texorpdfstring{$\GF$}{GF}: Far-Field Computation (P2G, G2P, and FFTs)}
\label{sec:sequential-far}

The far field is smooth and is defined in Fourier space:
\begin{subequations}\label{eq:GF}
\begin{equation} \label{eq:potentials-far}
  \u^\F(\x_i) = \sum_{j=1}^{\Ns} \Fou^{-1}\{\Gh{}^\F \cdot \f(\y_j)\}(\x_i - \y_j), 
\end{equation}
\begin{equation} \label{eq:GFk}
\Gh{}^F(\kv) =
 \frac{8\pi}{\lVert \kv \rVert^4}
\left(
\Id \lVert \kv \rVert^2 - \kv \otimes \kv
\right)
\E^{-\lVert \kv \rVert^2/(2\xi)^2}
\left(
1 + \frac{\lVert \kv \rVert^2}{(2\xi)^2}
\right),
\end{equation}
\begin{equation} \label{eq:inverse-transform}
\Fou^{-1}\{g\}(\r) = \frac{1}{L_1(2\pi)^2}
\sum_{\kappa_1 \in \K} \int_{\mathbb{R}^2} g(\kv) \E^{\I (\kv \cdot \r)}\,
\mathrm{d} \kappa_2\, \mathrm{d} \kappa_3,
\end{equation}
\end{subequations}
where $\kv=(\kappa_1,\kappa_2,\kappa_3)$,  $\K = \{2\pi \alpha / L_1 : \alpha \in \mathbb{Z}\}$,
$\Fou^{-1}$ is an inverse Fourier transform (series in the periodic direction and continuous in the free directions) and
$\Gh{}^\F$ the Fourier transform of $\GF=\G-\GN$. 
The support of $\Gh{}^\F$ grows with $\xi$, hence $\G^\F$ has slower spectral decay as $\xi$ increases.
For $\kappa_1=0$, the integral in \cref{eq:inverse-transform} is singular and is computed using the method by Vico et al.\ \cite{Vico_2016}. We omit details here and refer to \cite[Section~3.1]{Bagge_Tornberg_2023}.

A direct solve of \Cref{eq:GF} has $\mathcal{O}(N^2)$ complexity as the sources and targets are non-uniformly distributed in $\B$, hence fast Fourier transforms cannot be used.
This can be accelerated with a nonuniform FFT \cite{Barnett_Magland_Af_Klinteberg_2019}: spread points on a regular grid, use FFTs to convolve with $\GF$ in spectral space, and interpolate back. 
Let us assume that the regular grid has grid spacing $h$. The far-field potential $\u^\F$ is then computed using the following five computational kernels:
\begin{enumerate}[label=(\textbf{\arabic*})]
\item
\textbf{Particle-to-grid (P2G):} Spread the source strengths $\f(\y_j)$ to the grid points $\mathcal{G}_h:=\{\g_\ell\}$ of the uniform grid by convolving with a compactly supported window function $w(\r)$:
\begin{equation}
\label{eq:p2g}
\boldsymbol{\phi}_h(\g_\ell) := \sum_{j=1}^{\Ns} \sum_{\p \in \Ps} w(\g_\ell - \y_j + \p) \f(\y_j),
\qquad \ell=1,\ldots,\Ng,
\end{equation}
where $\Ng$ is the number of grid points and $\Ps$ is as in \cref{eq:potentials}.
\item
\textbf{FFT:} $\{\widehat{\boldsymbol{\phi}}_h(\kv_\ell)\}$ = \texttt{FFT3D} $\{\boldsymbol{\phi}_h(\g_\ell)\}$,
where $\{\kv_\ell\}$ are the grid points of the uniform grid in the frequency domain.
\item 
\textbf{Convolution with $\GF$ (CNV):}
Scale the values in Fourier space by the function $\Gh{}^\F$ and the Fourier transform $\widehat{w}$ of the window function:
\begin{equation}
\label{eq:sca}
    \widehat{\boldsymbol{v}}_h(\kv_\ell)
    := \frac{\Gh{}^\F(\kv_\ell)}{[\widehat{w}(\kv_\ell)]^2} \widehat{\boldsymbol{\phi}}_h(\kv_\ell).
\end{equation}
\item 
  \textbf{Inverse FFT:}  $\{\boldsymbol{v}_h(\g_\ell)\}$ = \texttt{IFFT3D} $\{\widehat{\boldsymbol{v}}_h(\kv_\ell)\}$.
\item 
  \textbf{Grid-to-particle (G2P):} Interpolate back to target points:
\begin{equation}
\label{eq:g2p}
    \u_h^\F(\x_i) := h^3 \sum_{\ell=1}^{\Ng} 
    \sum_{\p \in \Ps} w(\x_i - \g_\ell + \p)
    \boldsymbol{v}_h(\g_\ell).
\end{equation}
\end{enumerate}
The full potential is then approximated by $\u_h(\x_i) := \u^N(\x_i) + \u_h^\F(\x_i)$.
In our discussion, we refer to the combination of the FFT, CNV, and IFFT steps as the \textbf{Fourier grid convolution (FGC)} step.

The window function in the P2G and G2P steps is  $w(\r) = w_0(r_1)w_0(r_2)w_0(r_3)$, where $w_0$ is the truncated Kaiser--Bessel (KB)  function 
\begin{equation}
	w_0(r) = \begin{cases}
		\nicefrac{I_0\left(\beta\sqrt{1 - \nicefrac{r^2}{a_w^2}}\right)}{I_0(\beta)} & \text{if } r<a_w, \\ 
		0 & \text{otherwise},
	\end{cases}
\end{equation}
and $I_0$ is the modified Bessel function.
Since $I_0$ is expensive to evaluate (see \cref{tab:cycles}), $w_0$ approximated by a polynomial interpolation of degree $\nu$.
If $P$ denotes the number of grid points in the support of $w_0$,  then we set $P=2\nu+1$, $\beta = 2.5P$, $a_w = P h/2$. 
The truncated KB function was first used for the spectral Ewald method in \cite{Shamshirgar_2021} as an alternative to the classic Gaussian window function \cite{Lindbo_Tornberg_2011} with a smaller support needed for the same accuracy.
It has recently been shown that window functions based of direct polynomial approximation of prolates \cite[Appendix~A.6]{Jiang_Greengard_2025} has slightly smaller support than the truncated KB window; we leave such an implementation to future work. 
Computing \cref{eq:potentials-far} also requires extending the box $\B$ and the uniform grid in the free directions to approximate the free space conditions.
This is technical and not essential to our discussion; see~\cite[Section~4.5]{Bagge_Tornberg_2023} for details and a justification for the choices of $P$, $\beta$ and $a_w$. Since $\beta$ and $a_w$ are tied to $P$, the latter is the main window parameter, and it is set according to the accuracy requirements.

We summarize the entire spectral Ewald algorithm in \cref{alg:se}.
In blue is the complexity of the different steps, with $s$ being the \emph{average number of points per cell}.
 $\Ng$ is the total number of grid points in the uniform grid, $\Ng \sim h^{-3}$; $P$ is the window parameter related to the P2G and G2P interactions.

\begin{algorithm}
    \caption{Spectral Ewald method and its complexity}\label{alg:se}
    \begin{algorithmic}[1]
        \Require  $\B = [0, L_1) \times [0, L_2) \times [0, L_3)$; $\Ns$ source positions $\y_j \in \B$ and densities $\f(\y_j)$; $\Nt$ target positions $\x_i \in \B$;
        user-speified error tolerance $\varepsilon$.
        \Ensure Potential approximation $\u_h(\x_i)$ at target positions
        \State \textbf{(P2P)} Evaluate $\u^\N(\x_i)$ using \cref{eq:potentials-near}    {\color{blue} \hfill $ \mathcal{O}(3^d s \Nt)$}
        \State \textbf{(P2G)} $\f(\y_j) \rightarrow \boldsymbol{\phi}_h(\g_\ell)$ using \cref{eq:p2g} {\color{blue} \hfill $\mathcal{O}(\Ns P^d)$}
        \State \textbf{(FFT)}  $\{\boldsymbol{\phi}_h(\g_\ell)\} \rightarrow \{\widehat{\boldsymbol{\phi}}_h(\kv_\ell)\}$ {\color{blue} \hfill $\mathcal{O}(\Ng \log \Ng)$}
        \State \textbf{(CNV)}  $\widehat{\boldsymbol{\phi}}_h(\kv_\ell) \rightarrow \widehat{\boldsymbol{v}}_h(\kv_\ell)$ using \cref{eq:sca} {\color{blue} \hfill $\mathcal{O}(\Ng)$}
        \State \textbf{(IFFT)} $\{\widehat{\boldsymbol{v}}_h(\kv_\ell)\} \rightarrow \{\boldsymbol{v}_h(\g_\ell)\}$ {\color{blue} \hfill $\mathcal{O}(\Ng \log \Ng)$}
        \State \textbf{(G2P)} $\{\boldsymbol{v}_h(\g_\ell) \rightarrow \u_h^\F(\x_i)\}$ using \cref{eq:g2p} {\color{blue} \hfill $\mathcal{O}(\Nt P^d)$ }
        \State \textbf{(SUM)} Set $\u_h(\x_i) = \u^\N(\x_i) + \u_h^\F(\x_i)$ {\color{blue} \hfill $\mathcal{O}(\Nt)$ }
    \end{algorithmic}
\end{algorithm}
In the P2P kernel, $27$ is the number of  \emph{colleagues} of a cell --- the $3^d$ cells that share a face, edge, or point --- 
so each target interacts with $27s$ sources in its near field.
The $P^3$ factor in P2G is the cost of spreading one source density to the grid points in the support of the window function.

\subsection{Combined Single and Double Layers}\label{sec:doublelayer}
In many suspension problems with deformable particles, the Stokes potential is generated not only by Stokeslet sources (single layer densities), but also by stresslet sources (double layer densities).
That is, every source point then has an additional two vectors $\q(\y_j)$ and $\n(\y_j)$, and an additional term
$
\sum_{j=1}^{\Ns} \sum_{\p \in \Ps} \T(\x_i - \y_j + \p)  \big(
\q(\y_j) \otimes \n(\y_j)
\big)
$
gets added to the right-hand side of \cref{eq:potentials}.  The stresslet $\T$ is a symmetric 3-by-3-by-3 tensor
and is split again as $\T = \T^\N + \T^\F$, where the precise form of $\T^\N$ and $\T^\F$ is in \cite{Bagge_Tornberg_2023}.
If we combine the $\T$ and $\G$ kernels in one single call,
the amount of data throughout the algorithm changes as follows:
\begin{itemize}
\item The P2G step has a total of 12 components as input ($\n$, $\q$, $\y$, and $\f$ in $\mathbb{R}^3$), and 12 components as output ($\f$, $\q\otimes\n$).
\item The FFT has 12 components as input and output.
\item The CNV steps has 12 components as input, and 3 components as output (the $\G^\F+\T^\F$ combined $\widehat{\boldsymbol{v}}$).
\item The remaining steps (IFFT, G2P) are identical.
\end{itemize}
\emph{In the remaining discussion, all analysis is for this simultaneous treatment of Stokeslet and stresslet sources.}

\subsection{Accuracy and Parameters}
\label{sec:sequential-accuracy}
Here and in the remainder of the paper, without loss of generality, let $\B$ to be the unit cube so that $L_1=L_2=L_3=1$
and, for notational simplicity, let $\sum_j \| \f(\y_j) \|^2 = 1$.
We assume a nearly uniform distribution of points so that the number of points per cell $s$ does not vary much across cells.
Before proceeding further, we summarize the key algorithmic parameters:
\begin{itemize}
\item $\rc$ is the near-field cutoff radius, which we use as P2P cell spacing;
\item $s$ is the average number of points per cell;
\item $\xi$ controls near--far field decomposition scale and is related to $\rc$;
\item $h$ is the FFT uniform grid spacing, which controls accuracy and cost of $\Gh{}^\F$;
\item $P$ is the number of FFT grid points supported by the window function $w_0$ and controls accuracy and  cost of the P2G and G2P operations; 
\item $\eSE$ is an user-specified absolute error tolerance for the total approximation $\u_h$.
\end{itemize}

In~\cite{Bagge_Tornberg_2023} the root mean square error $E_{\mathrm{rms}} \coloneqq \sqrt{\nicefrac{1}{N_t}\sum_{i=1}^{N_t} |\u_h(\x_i) - \u(\x_i)|^2}$ is shown to be estimated by
\begin{equation}
\label{eq:ese}
	E_\mathrm{rms} \approx \underbrace{\sqrt{\frac{112}{9} \xi^4 \rc^3} \E^{-\xi^2\rc^2}}_{E^\text{N}}
    +
    \underbrace{\frac{4}{3 h} \sqrt{\frac{7}{2}} \E^{-\pi^2/(2 h \xi)^2}}
    _{E^\text{F}}
    +
    \underbrace{C(\xi) \E^{-2.5 P}}_{E^w}. 
\end{equation}
{\em The first term} $E^N$ is the near-field error due to the cutoff radius $\rc$.
{\em The second term} $E^F$ is the FFT error due to quadrature and truncation errors in the infinite series and integral in \cref{eq:GF}.
{\em The third term} $E^w$ is the error due to P2G and G2P interpolations.
We want to select all parameters to optimize runtime and ensure $E_\mathrm{rms}<\eSE$.%
\footnote{Besides the parameters we discuss here, there are two additional parameters related to the truncation of the integral in \cref{eq:inverse-transform}, but they can be also related to $\rc$. The overall methodology doesn't change.}

For a given $\rc$, we select $\xi:=\xi(\rc)$ so that $E^\text{N} < \eSE/3$.
Given $\xi$, we select $h:=h(\rc)$ and $P:=P(\rc)$ so that both $E^\text{F}$ and $E^w < \eSE/3$, respectively. 
Let $\mathcal{C}^N$ be the cell list set. 
Once we have these parameters, we set $|\mathcal{C}^N|=1/\rc^3$, i.e., the volume of $\B$ over the volume of a cell.
Thus, the average number of points per cell $s = \Ns/ |\mathcal{C}^N| = \Ns \rc^3$.
Similarly, the number of FFT grid points $\Ng = 1/h^3$.

Let $T(\Ng,\Ns, P,s)$ be the sum of the time complexity of all terms in \cref{alg:se}.
Since $\xi, P, h, \Ng,\Ns, P, s$ are implicit functions of $\rc$, both $T$ and $E_\mathrm{rms}$ are also implicit functions of $\rc$.
Thus, we can find $\rc$ by solving $\min_{\rc} T(\rc)$ subject to $E_\mathrm{rms}(\rc)\leq \eSE$.
Instead of the $T$ from \cref{alg:se} we can also use the GPU performance models with precise constants from \cref{sec:perf-mod}.
\emph{In conclusion, the algorithm has only one free parameter: the error tolerance $\eSE$, which is application dependent and set by the user}.%

%% file: sections/methods_gpu.tex
\section{GPU Acceleration}
\label{sec:acceleration}
We now discuss the design and implementations of GPU algorithms for the Ewald sum.
\Crefrange{sec:gpu-p2p}{sec:gpu-g2p} present GPU implementations of the \textbf{P2P}, \textbf{P2G}, and \textbf{G2P} steps. \Cref{sec:gpu-fgc} briefly addresses the \textbf{FFT}{}, \textbf{CNV} and the \textbf{IFFT} whose parallelization is more standard.
\Cref{sec:perf-port} discusses Kokkos-based performance portability to CPUs, 
\cref{sec:special} investigates the performance of different operations and special function evaluations on GPUs,
and \cref{sec:perf-mod} presents performance models for each stage of the Spectral Ewald method.

\begin{figure}
    \centering
	\includegraphics[width=0.6\textwidth]{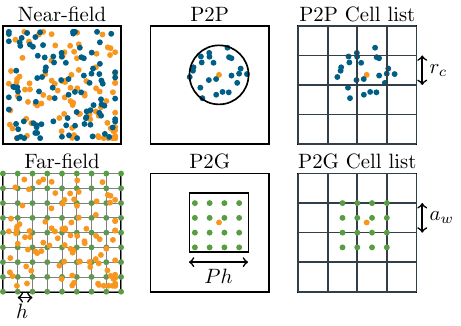}
    \caption{Local structure for near (\emph{top row}) and far (\emph{bottom row}) interactions.
	    \emph{Left column:} source particles $\y_j$ in orange (\emph{top/bottom}), target particles $\x_i$ in blue (\emph{top}) and the regular grid points $\g_\ell$ in green (\emph{bottom}), with grid spacing $h$.
	    \emph{Middle column:} interaction radius for a source particle $\y_j$.
	    The P2P interaction radius is $r_c$.
	    The P2G interaction radius is $a_w = Ph/2$, where $P$ is the number of grid points needed to support $w_0(r)$.
	    \emph{Right column:} cell list partition; the interaction radius requires each source cell to interact with each $3^d$ of its colleagues. 
	    Each source cell is mapped to a GPU block.
	    The G2P kernel is similar to P2G; the sources $\y_j$ are replaced with targets $\x_i$.
    }
    \label{fig:cell-list}
    \end{figure}

As mentioned in \cref{sec:sequential}, \textbf{P2P} requires a cell list $\CN$ for efficient computation. To promote \emph{coalesced} memory access on the GPUs, the cell list is \emph{padded} so that each cell contains the same number of points. In addition, padded cell lists can accelerate \emph{each} of the \textbf{P2P}, \textbf{P2G}, and \textbf{G2P} steps by assigning a GPU block $B$ to a cell $\beta$, ensuring that each thread $(t_x,t_y) \in B$ processes spatially close data. \Cref{fig:cell-list} summarizes the cell lists for \textbf{P2P}, \textbf{P2G}, and \textbf{G2P}.

The construction of a padded cell list can be implemented in two passes through an unsorted particle array. The first pass, carried out in parallel, counts the maximum number of particles per cell $s$ and the number of nonempty cells $N_c$. Then the cell-list array of shape $(s\, N_c, d)$ is allocated --- where $d$ is the dimension of the data to be stored --- as well as a cell counter of size $N_c$. The second pass reads the particle data from the unsorted array to the cell-list array by assigning one thread to one particle and logging the cell-list positions with \Call{AtomicFetchAdd}{} updates to the cell counter. Cells are ordered via a linear index of a 3D partition which enables efficient indexing and memory access on GPU.  In the following, we refer to the 1D GPU block size as $b_x$ and the 2D GPU block size as $(b_x, b_y)$.
\subsection{Particle-to-Particle Interaction (P2P)}
\label{sec:gpu-p2p}

Near-field interactions are a standard component of parallel $N$-body libraries. 
Typically, a target cell $\beta \in \CN$ is assigned to a GPU block; each target point $\x_m \in \beta$ is assigned to a thread, and a serial loop iterates over the sources points $\y_j$ in the colleague cells $\alpha$ of $\beta$.
We call this approach \Call{P2P-GM-1D}{}: a 1D parallel loop over targets, with sources read from global memory (GM).
We can also store sources in shared memory (SM) as each source $\y_j \in \alpha$ is shared by every target $\x_i \in \beta$.
To express this operation, we let $\Call{LoadShmem}$ be a function that loads its arguments from the GPU memory to the shared memory.
Another algorithmic variant uses a 2D thread grid and parallelizes on both target and sources. 
Assuming that GPU threads are indexed by $(t_x,t_y)$,  $t_x$ is assigned to a target $x_i$ and
each $t_y$ handles \emph{chunks} of source points in the cell $\alpha$.
We refer to these chunks as $\alpha_c$.
Thread $(t_x,t_y)$ accumulates all the interactions between $\x_i$  and $\y_j \in \alpha_c$ to a \emph{chunk potential} $\u_c(\x_i)$.
Finally, the threads in the row $t_x$ perform a reduction to compute $\u^\N(\x_i) = \sum_c u_c(\x_i)$.

As in the 1D case, sources can be read from global or shared memory.
The 2D shared-memory strategy offers increased parallelism over the sources at the cost of synchronization and diminished target-wise parallelism,
as the total number of threads is fixed but more threads are allocated per target. 
All four variants follow the template in \cref{alg:p2p-sm-2d}, which is the 2D shared memory variant.
{\Call{P2P-GM-1D}{}} doesn't use lines 4--5; $\alpha_c = \alpha$ so there is no loop in line 3; and there is no reduction.  
{\Call{P2P-SM-1D}{}} uses lines 4--5; $\alpha_c = \alpha$ so there is no loop on line 3; and there is no reduction.
{\Call{P2P-GM-2D}{}} doesn't uses lines 4--5.
{\Call{P2P-SM-2D}{}} is shown below.

\begin{algorithm}
  \footnotesize
    \caption{P2P-SM-2D}\label{alg:p2p-sm-2d}
    \begin{algorithmic}[1]\setcounter{ALG@line}{-1}
        \Require $\CN$
        \Ensure $\u_h^N$
        \State Assign target cell $\beta$ to a GPU block. 
        \ParFor{$\x_i \in \beta$} \Comment{Assign $\x_i$ to thread $t_x$.} 
            \ParFor{$c \gets 1, \dots, \lceil s/|\alpha_c| \rceil$} \Comment{Assign $c$ to thread $t_y$}
            \State $\u_c \gets 0$
            \State $\forall \y_j \in \alpha_c, \Call{LoadShmem}{\y_j, \f(\y_j), \q(\y_j), \n(\y_j)}$
            \State \Call{\_\_synchthreads}{ }
                \For{$\y_j \in \alpha_c$}
                    \State $\u_{\mathrm{c}} \gets \u_{\mathrm{c}} + \G^N(\x_i-\y_j+\p)  \f(\y_j)$
                    \State $\u_{\mathrm{c}} \gets \u_{\mathrm{c}} + T^N(\x_i - \y_j+\p)  (\q(\y_j)\n(\y_j))$
                \EndFor
                \State \Call{\_\_synchthreads}{ }
            \EndParFor
            \State{$\u_h^N(\x_i) \gets  \Call{RowReduce} {\u_{\mathrm{c}}}$} \Comment{row-wise reduction}
        \EndParFor
    \end{algorithmic}
\end{algorithm}

\subsection{Particle-to-Grid Interaction (P2G)}\label{sec:gpu-p2g}
\begin{figure}
    \centering
	\includegraphics[width=0.6\textwidth]{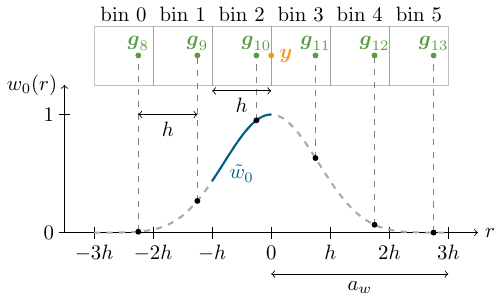}
    \caption{Evaluation of the one-dimensional piecewise approximation $\widetilde{w}_0$ (solid curve) to the KB window function $w_0$ (dashed curve) at a source $y$ on the Fourier grid. 
      An anchor point $a =  g_8$ is found and $\tilde w_0$ is evaluated at bin 2 (i.e., $g_{8+2}$) 
      $\widetilde{w}_0(\g_{8+2} - y)=\widetilde{w}_0(2h - (y-a))$. (Here, $P=6$.)}
    \label{fig:window}
\end{figure}

The P2G interaction interpolates the density from the source points to the grid.
Its implementation involves several design decisions:
Do we loop over grid points gathering contributions or do we loop over source points scattering contributions?
Do we incur the costs associated with reordering source points in order to optimize streaming memory access or not?
Can we avoid conditionals to find the correct interpolant?
To explore these questions, we consider four algorithmic variants. 
Before describing them, it suffices to consider the P2G steps in 1D as the steps are identical across dimensions.
Let $y$ be a source point.
We use $\widetilde{w}_0$ to denote the polynomial approximation of $w_0$. 
To compute $\widetilde{w}_0(r)$ we need to find the $r$-axis interval $y$ belongs to; we use the term \emph{bin} to refer to such $r$-intervals. 
\cref{fig:window} summarizes the calculation that comprises the following steps:
\begin{enumerate*}[label=(\textbf{\arabic*})]
    \item find a \emph{bin anchor point} $a \gets \lfloor y/h \rfloor - (P/2 - 1)$ on the grid;
    \item compute a local distance $\delta = a - y$;
    \item compute $\widetilde{w}_0(a+lh - y)$ for each \emph{bin} $l \gets 0, \dots, P-1$.
\end{enumerate*}

Two of the variants require reordering the source points. To do so, we create a new cell list \CF.
It is constructed identically to \CN, but with a cell width equal to $a_w=Ph/2$.
There is no need to explicitly reorder the grid points as their positions are computed analytically.
Writing the grid values for the density does require uncoalesced writes; but the writes enjoy some locality since source points share grid points.
Also, let $\Cell$ denote the user-defined input order, which in general can be considered random; and let $\mathcal{G}$ indicate the Fourier grid parameters and array.

Then, the four variants, all using 1D GPU blocks, are as follows:
\begin{itemize}[topsep=0pt,itemsep=-1ex,partopsep=1ex,parsep=1ex]

\item {\bf Baseline} (\Call{P2G-base}{}, \cref{alg:p2g-source}): Loop over source points in \Cell\ order; compute contribution to grid points and use atomic add to update $\boldsymbol{\phi}$.

\item {\bf Source-ordered} (\Call{P2G-source}{}, \cref{alg:p2g-source}): Same as above, but source points are read in \CF\ order.

\item {\bf Grid-ordered} (\Call{P2G-grid}{}, \cref{alg:p2g-grid}): Loop over \emph{grid points} to  avoid atomics; requires repeated reads of source points, hence repeated evaluations of the window function for a source. 

\item {\bf Grid-hybrid} (\Call{P2G-hybrid}{}, \cref{alg:p2g-hybrid}): Use shared memory for the source points and traverse all neighbors of a source point to spread the density.
\end{itemize}

\begin{algorithm}
\footnotesize    
    \caption{P2G-base / P2G-source}\label{alg:p2g-source}
    \begin{algorithmic}[1]\setcounter{ALG@line}{-1}
        \Require $\Cell/ \Cell^F$ (source ordering), $\mathcal{G}$ (grid)
        \Ensure $\boldsymbol{\phi}(\mathcal{G})$
        \State Assign  $\alpha$ in $\mathcal{C}/\mathcal{C^\F }$ to a GPU block.   \Comment{Base / Source}
        \ParFor{$\y_j \in \alpha$} \Comment{Assign $\y_j$ to a thread}
            \State $ \Call{LoadGlb2Reg}{\y_j, \f(\y_j), \q(\y_j), \n(\y_j)}$
            \State $\boldsymbol{\sigma} \gets [\f(\y_j), \mathrm{vec}(\q(\y_j) \otimes\n(\y_j)]^T$
            \State $\a \gets \lfloor\y_j/h\rfloor - (P/2)  + 1$ \Comment{Compute grid anchor}
            \State $\boldsymbol{\delta} \gets \y_j/h - \a - \p$ \Comment{Compute distance offset}
            \For{$l \gets 0, \dots, P-1$} \Comment{Loop over bin $l$}
                \State $w_1[l] \gets \widetilde{w}_0(lh-\delta_1)$
                \State $w_2[l] \gets \widetilde{w}_0(lh-\delta_2)$
                \State $w_3[l] \gets \widetilde{w}_0(lh-\delta_3)$
            \EndFor
            \For{$k \gets 0, \dots, P-1$}
                \For{$l \gets 0, \dots, P-1$}
                    \For{$m \gets 0, \dots, P-1$}
                        \State $\g_m \gets \a + [k,l,m]^T $
                        \State $w \gets w_1[k]\, w_2[l] \, w_3[m]$ 
                        \State $\Call{AtomicAdd}{\boldsymbol{\phi}(g_i), w \, \boldsymbol{\sigma}}$ \Comment{$\mathcal{O}(P^3N_s)$ atomics}
                    \EndFor
                \EndFor
            \EndFor
        \EndParFor
    \end{algorithmic}
\end{algorithm}

\Call{P2G-base}{}  is easy to implement, but
it incurs $\mathcal{O}(P^3\Ns)$ atomic adds and the source point reads are completely random, although it is possible that this randomness reduces the overhead from atomics.  \Call{P2G-source}{} attempts to improve the source reads, but suffers from similar problems. Both versions avoid thread divergence, since threads in a warp evaluate $\widetilde{w}_0$ in the same bin (lines 6–9).

In contrast, \Call{P2G-grid}{}, which parallelizes over the grid points $\g_\ell$ and serially loops over the sources $\y_j$, completely avoids atomics. 
However, due to the piecewise-polynomial approximation of $\widetilde{w}_0$, thread divergence is introduced into the serial loop; for a fixed $\g_\ell$, the bin $[k,l,m] \gets \g_\ell - \a$, where $\a$ is an anchor point dependent on $\y_j$, will vary randomly for given sources $\y_j$, conflicting with the GPU execution model. 
Additionally, parallelizing over the grid points increases the work of evaluating $w$ at a source $\y_j$ from $\mathcal{O}(2P^3)$ to $\mathcal{O}(3f(\widetilde{w}_0)P^3)$, where $f(\widetilde{w}_0)$ is the cost of computing $\widetilde{w}_0$.
Finally, we load the same source repeatedly --- in the worst case $\mathcal{O}(3^d)$ times --- ignoring caching.

\begin{algorithm}
\footnotesize  
    \caption{P2G-grid}\label{alg:p2g-grid}
    \begin{algorithmic}[1]\setcounter{ALG@line}{-1}
        \Require $\CF$, $\mathcal{G}$
        \Ensure $\boldsymbol{\phi}(\mathcal{G})$
        \State Assign target box $\alpha$ to a GPU block.
        \ParFor{$\g_\ell$} \Comment{Assign $\g_\ell$ to a thread}
            \State $\boldsymbol{\phi}_{\mathrm{loc}}\gets 0$
            \For{$\beta \in $ {\tt Colleagues}$(\alpha)$}
                \For{$\y_j \in \alpha$}
                    \State $ \Call{LoadGlobal}{\y_j, \f(\y_j), \q(\y_j), \n(\y_j)}$
                    \State $\boldsymbol{\sigma} \gets [\f(\y_j), \mathrm{vec}(\q(\y_j) \otimes\n(\y_j)]^T$
                    \State $\a \gets \lfloor\y_j/h\rfloor - (P/2) + 1$
                    \State $\boldsymbol{\delta} \gets \a - \y_j/h - \p$
                    \State $[k,l,m]^T \gets \g_\ell - \a$
                    \State $w \gets \widetilde{w}_0(kh-\delta_1)\,\widetilde{w}_0(lh-\delta_2)\,\widetilde{w}_0(mh-\delta_3) $ \Comment{Thread divergence!}
                    \State $\boldsymbol{\phi}_{\mathrm{loc}}\gets \boldsymbol{\phi}_{\mathrm{loc}} + w  \boldsymbol{\sigma}$
                \EndFor
            \EndFor
            \State $\boldsymbol{\phi}(\g_\ell)\gets \boldsymbol{\phi}_{\mathrm{loc}}$
        \EndParFor
    \end{algorithmic}
\end{algorithm}

The drawbacks of the first three variants motivate a two-stage \emph{hybrid} parallelization scheme.
Like in the \Call{P2P-base}{} and \Call{P2G-source}{} methods, each source cell $\alpha$ is assigned to a GPU block (or more if it has too many source points).
Then for each source, we load the densities to the shared memory and save them for the loop beginning at line 15.
Lines 4--11 also store the interpolation coefficients per point in shared memory.
Once we synchronize on line 14, we assign each thread an index $t_x$ between $0$ and $(P/2)^3$ using the thread id. 
(If $(P/2)^3 < b_x$, we have some idle threads;  if $(P/2)^3>b_x$, we stride by $b_x$ to assign multiple indices to a thread.)
In line 16, we begin a sequential loop over the colleague cells. 
Regular grid ordering guarantees that every colleague cell $\beta$ has $(P/2)^3$ grid points. 
As such, we define a \emph{local lexicographical ordering} on grid points $\g_\ell \in \beta$, and use integer arithmetic to determine the global grid coordinates of $\g_\ell$ from the thread index $t_x$ in a function we call \Call{Index2Grid}{} to assign each thread a grid point $\g_\ell$. 
Lines 18--26 mirror the computation in lines 4--12 of \Call{P2G-grid}{}, but now $\a_j, f(\y_j), \q(\y_j), \n(\y_j)$ are read from \emph{shared memory} as well as the 1D window values $w_1[k], w_2[l], w_3[m]$, avoiding the thread divergence in \Call{P2G-grid}{} line 10. 
On line 27, we \emph{atomically} write $\boldsymbol{\phi}(\g_\ell)$, as $\mathcal{O}(3^d)$ source cells $\alpha$ share a given colleague cell.
\emph{In summary, \Call{P2G-hybrid}{}\relax\ minimizes atomics,  reuses source points, and avoids repeated computations of the window function.
In fact, we experimented replacing line 27 with a different variant that first writes at different locations and then calls an additional kernel for reduction, but it was slower.}

\begin{algorithm}
\footnotesize    
    \caption{P2G-hybrid}\label{alg:p2g-hybrid}
    \begin{algorithmic}[1]\setcounter{ALG@line}{-1}
        \Require $\CF$, $\mathcal{G}$
        \Ensure $\boldsymbol{\phi}(\mathcal{G})$
        \State Assign source box $\alpha$ to a GPU block
        \ParFor{$\y_j \in \alpha$} \Comment{Assign $\y_j$ to a thread}
            \State $\Call{LoadGlobal}{\y_j}$
            \State $ \Call{StoreGlobal2Shmem}{\f(\y_j), \q(\y_j), \n(\y_j)}$
            \State $\a_j \gets \lfloor\y_j/h\rfloor - (P/2)  + 1$
            \State $\Call{StoreShmem}{\a_j}$
            \State $\boldsymbol{\delta} \gets \y_j/h - \a_j - \p$ 
            \For{$l \gets 0, \dots, P-1$} 
                \State $w_1[l] \gets \widetilde{w}_0(lh-\delta_1)$
                \State $w_2[l] \gets \widetilde{w}_0(lh-\delta_2)$
                \State $w_3[l]\gets \widetilde{w}_0(lh-\delta_3)$
                \State $\Call{StoreShmem}{w_1[l], w_2[l], w_3[l]}$
            \EndFor
        \EndParFor
        \State \Call{\_\_synchthreads}{ }
        \ParFor{$t_x \gets 0, \dots, (P/2)^3 - 1$} \Comment{Use $(P/2)^3$ threads}
            \For{$\beta \in $ {\tt Colleagues}$(\alpha)$}
                \State $\g_\ell \gets \Call{Index2Grid}{t_x, \alpha}$ \Comment{Assign $\g_\ell$ to a thread}
                \State $\boldsymbol{\phi}_{\mathrm{loc}}\gets 0$
                \For{$\y_j \in \alpha$} 
                    \State $ \Call{LoadShmem}{\a_j, \f(\y_j), \q(\y_j), \n(\y_j)}$
                    \State $\boldsymbol{\sigma} \gets [\f(\y_j), \mathrm{vec}(\q(\y_j) \otimes\n(\y_j)]^T$
                    \State $[k,l,m]^T \gets \g_\ell - \a_j$
                    \State $\Call{LoadShmem}{w_1[k], w_2[l], w_3[m]}$
                    \State $w \gets w_1[k]\, w_2[l]\, w_3[m]$
                    \State $\boldsymbol{\phi}_{\mathrm{loc}}\gets \boldsymbol{\phi}_{\mathrm{loc}} + w \, \boldsymbol{\sigma}$
                \EndFor
                \State $\Call{AtomicAdd}{\boldsymbol{\phi}(\g_\ell), \boldsymbol{\phi}_{\mathrm{loc}}}$ \Comment{$\mathcal{O}(27N_s)$ atomics}
            \EndFor
        \EndParFor
    \end{algorithmic}
\end{algorithm}

\subsection{Grid-to-Particle Interaction (G2P)} \label{sec:gpu-g2p}
This operation can be thought of as the ``\emph{transpose}'' of the P2G interaction.
Here we interpolate the potential from the grid points to the target points. 
For brevity, we only outline this interaction.
We consider two variants {\Call{G2P-base}{}} and {\Call{G2P-target}{}}.
The first variant is essentially identical to \Call{P2G-base}{}, the difference being that in line 1 in \cref{alg:p2g-source} we parallel loop over \emph{target points}.
Then each thread sequentially loops over the grid points in the support of KB window centered at the target, and interpolates the potential at the target point.
Unlike the P2G kernel, no atomics are needed as each thread adds to each target point potential.
As such, in \Call{G2P-target}{} the only optimization we make is to reorder the targets using $\CF$ to improve data locality.

\subsection{Fourier Grid Convolution (FGC)}
\label{sec:gpu-fgc}
We now discuss the computations in Fourier space. Single GPU FFT and IFFT are done using vendor libraries.
The convolution with  $\Gh{}^\F$ is diagonal; we assign a thread to each frequency $\boldsymbol{\kappa}_\ell$ and compute the multiplication in \cref{eq:sca}.
Moreover, since each of $\boldsymbol{\phi}$, $\G^\F$, and $w$ is real-valued, their Fourier transforms exhibit conjugate symmetry.
As a result, the convolution only needs to be computed for the \emph{nonnegative} wavenumbers $\kappa_{\ell}$, effectively reducing the computational workload by half compared to the general nonsymmetric case.

\subsection{Portable Implementation Using Kokkos Abstractions}
\label{sec:perf-port}
All interactions and their variants have been implemented using Kokkos abstractions via the Python interface PyKokkos.
We map each GPU block to a Kokkos team and utilize Kokkos scratchpad memory to emulate managed shared memory.
For 2D kernels, thread rows are mapped to Kokkos vector ranges, which are guaranteed to execute within the same warp.
Kokkos also provides built-in support for vector reductions, which we leverage to implement \Call{RowReduce}{} in \cref{alg:p2p-sm-2d}.
This level of GPU abstraction enables seamless portability of our kernels across both CUDA, HIP, and CPU backends.
A shortcoming for CPUs is that we do not use explicit vectorization but rather rely on the compiler, although Kokkos does offer SIMD data types for abstracting vector intrinsics. 
When targeting the CPU, Kokkos teams (originally corresponding to GPU blocks) are mapped to OpenMP threads.
Accordingly, we set the block size to one, effectively restructuring the cell list so that each team operates on a single particle per cell.

%% file: sections/methods_erf.tex
\section{Special Function Evaluation}
\label{sec:special}

The scheme involves a number of special mathematical functions such as $\erfc(x)$, $\E^{-x}$ and $\sinh(x)$. 
To estimate the cost of these functions, we provide \cref{tab:cycles}, evaluating the double-precision CUDA functions with random input numbers in the range $[10^{-3}, 7]$, which is typical for our method, on an NVIDIA H200 GPU using NVIDIA Nsight Compute. 
Performance is measured in the average number of millicycles per operation and normalized by the half the average number of millicycles per fused multiply-add (DFMA) to obtain a \emph{FLOP-equivalent} metric. 
For example, a division (DIV) requires on average 16.6 as many cycles as a multiplication (DMUL), and is computed using 15 double precision (and 2 single precision) flops. For some functions (e.g., $J_1$), the cycle count is much larger than expected from the flop count, likely due to thread divergence.

\begin{table}[t]
\centering
\caption{Observed millicycles per operation for a selection of FP64 CUDA functions on the NVIDIA H200 GPU, measured with NVIDIA Nsight Compute. Flop-equivalents are computed by dividing the millicycles of each operation by those of a DFMA and multiplying by 2. Actual FP64 flops are obtained from Nsight Compute; additional FP32 flops in parenthesis if present. Shaded rows are machine instructions.}
\label{tab:cycles}
\setlength{\aboverulesep}{0pt}
\setlength{\belowrulesep}{0pt}
\setlength{\tabcolsep}{3pt}
\begin{tabular}{lSSc}
  \toprule
  \footnotesize
\textbf{Operation} & \textbf{Millicycles} & \textbf{Flop-equivalents} & \textbf{Actual flops} \\
\midrule
\rowcolor[gray]{0.9}DFMA & 0.121 & 2 & 2 \\
\rowcolor[gray]{0.9}DADD & 0.127 & 2.10 & 1 \\
\rowcolor[gray]{0.9}DMUL & 0.127 & 2.10 & 1 \\
DIV & 2.11 & 35.0 & 15 (2) \\
SQRT & 2.50 & 41.5 & 13 \\
RSQRT & 1.82 & 30.2 & 8 \\
EXP & 3.07 & 50.8 & 30 \\
ERF & 9.33 & 155 & 78 (1) \\
ERFC & 16.1 & 267 & 112 \\
SINH & 21.8 & 361 & 49 (2) \\
I0 & 5.00 & 82.9 & 36 \\
J0 & 15.8 & 262 & 36 \\
J1 & 39.4 & 654 & 36 \\
\bottomrule
\end{tabular}
\end{table}

As mentioned in \cref{sec:sequential-far}, the window function $w_0(r)$ is approximated by a piecewise polynomial of degree $\nu$ to avoid the cost of evaluating the Bessel function $I_0$. 
A separate polynomial is used in each bin ($r$-axis interval) of width $h$ shown in \cref{fig:window}, and these polynomials are precomputed, as described in \cite{Shamshirgar_2021}. 
The polynomial approximation of $w_0$ can be evaluated using Horner's rule, requiring only $\nu$ DFMA operations per evaluation point.

\begin{figure}
    \centering
\includegraphics[width=0.6\linewidth,trim=7mm 8mm 8mm 7mm,clip]{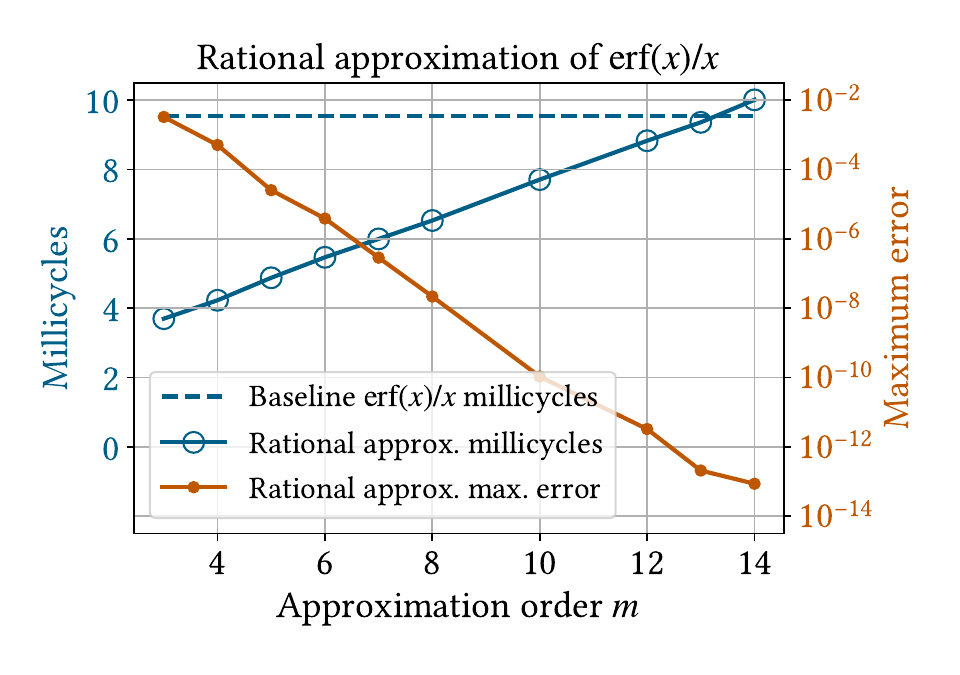}
    \caption{Millicycles per evaluation and maximum error for rational approximation of $\erf(x)/x$ on the interval $[0,6]$.}
    \label{fig:rational-runtime}
\end{figure}

The $\erfc(x)/x$ function that appears in the P2P step, cf.\ \cref{eq:g-near}, can also be approximated. First note that $\erfc(x)/x=1/x-\erf(x)/x$, and $\erf(x)/x$ is a smooth and slowly decaying function that is suitable for rational function approximation. We apply the AAA algorithm \cite{Nakatsukasa_2018,Nakatsukasa_2020} to construct a barycentric rational approximant $r(x) = \sum_{j=1}^m \frac{w_j f_j}{x-z_j} / \sum_{j=1}^m \frac{w_j}{x-z_j}$ of $\erf(x)/x$ for $x \in [0, 6]$. However, the barycentric form has $m+1$ divisions, which is far too expensive. 
Instead, we convert the output of the AAA algorithm to the polynomial form $r(x)=p(x)/q(x)$ where both $p$ and $q$ are polynomials of degree $m-1$, allowing us to use $2(m-1)$ DFMAs and a single division to evaluate $r(x)$. 
The millicycles observed and the maximum error as a function of $m$ are shown in \cref{fig:rational-runtime}. At $m=8$, the error is $10^{-8}$ and the speedup of using rational approximation instead of the CUDA function is 1.46$\times$.

Empirical measurements suggest that around 22.5\% of the P2P runtime is taken up by $\erf(x)/x$ evaluation, when using the CUDA $\erf$ function. The speedup for the entire P2P runtime when using the rational approximation with $m=8$ would then be 1.08$\times$.

%% file: sections/methods_model.tex
\section{Performance Models}
\label{sec:perf-mod}
We construct GPU performance models --- to evaluate the algorithms described in \cref{sec:acceleration} --- that account for global memory, L2 cache, registers, and shared (or L1) memory. 
The models are constructed with increasing levels of cache reuse and are denoted by 
\begin{itemize}
	\item $T_0$ (used to model P2G/G2P): A 0-L2 model with no L2 or L1 cache reuse but with managed shared memory and an assumption of coalesced high-bandwidth memory (HBM) access for bandwidth values. 
  \item $T_\infty$: Same memory characteristics as $T_0$ with an additional infinite L2 cache; data is loaded from HBM once and reused from L2 thereafter.
  \item $T_{\infty, l_1}$: Same memory characteristics as $T_\infty$ with an additional infinite L1 cache; data is reused from L1 once prefetched from L2.
  \item $T_f$ (used to model P2P): A compute-only model that accounts solely for floating-point operations, ignoring all data movement costs.
\end{itemize}
The models are defined by the following variables:
\begin{itemize}
	\item $\tau_f$, the machine time/FLOP; 
	\item $\tau_m$, the global memory time/byte (load or store); 
	\item $\tau^1_m$, the L1 memory time/byte (load or store); 
	\item $\tau^2_m$, the L2 memory time/byte (load or store); 
	\item $l_1 = \tau^1_m/\tau_m$ and $l_2 = \tau^2_m/\tau_m$; 
	\item $f$, the algorithm total flop count; 
	\item $m$, sum of bytes for all the input data; 
	\item $m_G$ and $m_S$, the memory load/stores from the global and shared memory, respectively.
\end{itemize}
Model specific memory trafic estimates are given by $m_0$, $m_\infty$, and $m_{\infty,l_1}$, defined as $m_0 = m_G + l_1 m_S$,  $m_{\infty} = m + l_2 m_G + l_1 m_S$ and $m_{\infty,l_1} = m + l_1 m_G + l_1 m_S$.
The models are explicitly expressed as $T_0 = \tau_f f + \tau_m m_0$,  $T_\infty = \tau_f f + \tau_m m_\infty$, $T_{\infty, l_1} = \tau_f f + \tau_m m_{\infty, l_1}$ and $T_f = \tau_f f$.
Arithmetic intensity is given by $q_k = \nicefrac{f}{m_k}$ for $k \in \{0,\ \infty,\ \infty,l_1\}$. If $q_k \geq \nicefrac{\tau_m}{\tau_f}$, the kernel is compute-bound under model $k$; otherwise, it is memory-bound.

\Cref{sec:perf-p2p,sec:perf-p2g,sec:perf-g2p} provide the analytical value of $f$, $m_0$, $m_\infty$ and $m_{\infty,l_1}$ for each kernel (P2P, P2G, and G2P) and compare model predictions with empirical GPU runtimes.
Sources and targets are independently  generated by a nearly uniform distribution in the unit cube and
all numerical quantities and arithmetic are in double precision, unless otherwise noted.

\begin{figure}
    \centering
\includegraphics[width=0.7\linewidth]{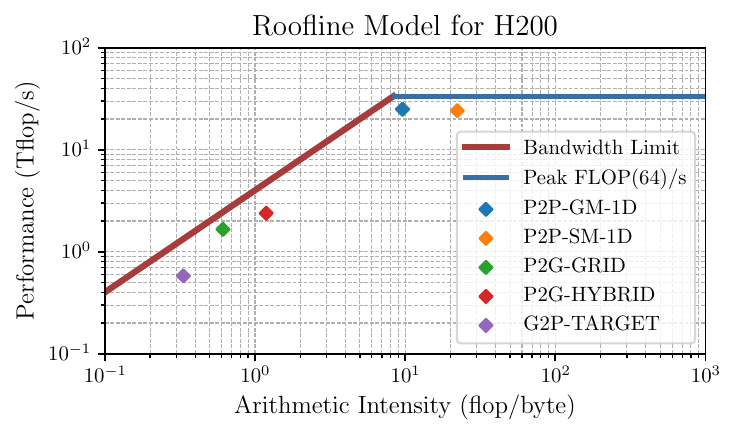}
\caption{Analytical roofline and flops of \cref{alg:se} on the NVIDIA H200 GPU for $N = 4\times10^6$, $s = 224$, $\epsilon = 10^{-9}$.
  $T_{\infty}$ model is used for P2P kernels,
  $T_0$ for P2G and G2P kernels.}
    \label{fig:roof}
\end{figure}

We provide a brief summary of the accuracy of our models:
\begin{itemize}
    \item {\bf P2P} (\cref{tab:p2p_mop_all}): The SM-1D kernel has measured times slightly lower than $T_{\infty, l_1}$, which we attribute to partial overlap between memory transfers and computation.
    \item {\bf P2G, G2P}  (\cref{tab:p2g_time_model}): Measured times exceed the $T_0$ model as the model assumes coalesced accesses to HBM, but both P2G and G2P have irregular global memory accesses.
\end{itemize}
Efficiency for the P2P kernel is calculated using $T_f$ as a baseline, while efficiency for the P2G and G2P kernels are
calculated using $T_0$ as a baseline.
The models are optimistic; achieving high efficiency compared to the baseline provides strong evidence that our kernels are highly optimized.
\cref{fig:roof} presents a roofline analysis of the different kernels.

More accurate models such as Volkov~\cite{volkov2016understanding} can be used, but our models are designed to provide a reasonable approximation to observed performance behavior by capturing the dominant data movement and compute costs in each kernel.

\subsection{P2P}\label{sec:perf-p2p}
The  P2P kernel's flop count is given by $f_\text{P2P} = 27N_ts\left(C_{\mathrm{p2p}}\nicefrac{\pi}{6} + 14\right)$,
where $s$ is the number of sources in one cell and $C_{\mathrm{p2p}}$ varies across machines due to differing costs for special function evaluations (see \cref{tab:cycles}).
For instance, on the NVIDIA A100 GPU, $C_{\mathrm{p2p}} = 318$, while on the NVIDIA H200 GPU, $C_{\mathrm{p2p}} = 414$.
The load operations for the P2P kernel include the target points coordinates  ($3N_t$ doubles) and the source points coordinates, normal vectors, and single-layer and double-layer densities ($12N_s$ doubles).
The outputs are the potentials at the target points ($3N_t$ doubles).
Thus, $m = 48N_t + 96N_s$ for all P2P variants, while $m_G$ and $m_S$ depend on the variant. 
In \Call{GM-1D}{}, each thread loads one target and $27s$ sources and globally stores one target.
Thus, $m_G^\mathrm{GM-1D} = N_t(48 + 96\cdot27s)$, and $m_S^\mathrm{GM-1D} = 0$.
In \Call{SM-1D}{}, given a GPU block of size $b_x \leq s$, each block moves $b_x$ targets and $27s$ sources (in chunks of size $b_x$) into shared memory.
Thus, $m_G^\mathrm{SM-1D} = N_t(48 + \nicefrac{96\cdot27s}{b_x})$ and $m_S^\mathrm{SM-1D} = N_t(24 + 96\cdot27s)$.
For the 2D variants, let $(b_x, b_y)$ denote the size of the 2D block.
Then $m_G^\mathrm{GM-2D} = N_t(24b_y + 24 + 96\cdot27s)$ and $m_S^\mathrm{GM-2D} = 0$; $m_G^\mathrm{SM-2D} = N_t(48 + \nicefrac{96\cdot27s}{b_x})$ and $m_S^\mathrm{SM-2D} = N_t(24b_y + 96\cdot27s)$.
\Cref{tab:p2p_mop_all} summarizes the memory operations.
\cref{tab:p2p_time_model} shows the empirical and analytical runtime under different performance models for different P2P algorithmic variants.
The results suggest that the $T_{\infty}$ model is somewhat pessimistic, likely because it does not account for prefetching from L2 to the system managed L1. 
This observation motivates the inclusion of the $T_{\infty, l_1}$ and $T_f$ models at the beginning of this section.
Notably, the SM-1D variant achieves measured times even lower than $T_{\infty, l_1}$, indicating that the kernel is highly optimized.

\begin{table}
\centering
\caption{Analytical memory operation (MOP) counts for P2P algorithmic variants under the $T_{\infty}$ model. GPU block size is denoted $b_x$ for a 1D block and $(b_x, b_y)$ for a 2D block. Here, $m=48N_t + 96N_s$.  The MOP value for the $T_0$ model ($m_0$) is obtained by setting $m = 0$ and $l_2 = 1$, while the MOP value for the $T_{\infty, l_1}$ model ($m_{\infty, l_1}$) is obtained by setting $l_2 = l_1$.}
\label{tab:p2p_mop_all}
\setlength{\aboverulesep}{0pt}
\setlength{\belowrulesep}{0pt}
\setlength{\tabcolsep}{3pt}
\begin{tabular}{lll}
\toprule
\multirow{4}{*}{$m_{\infty}$}
  & \Call{GM-1D}{} & $m + l_2 N_t(48 + 96 \cdot 27s)$ \\
  & \Call{GM-2D}{} & $m + l_2 N_t(24b_y + 24 + 96 \cdot 27s)$ \\
  & \Call{SM-1D}{} & $m + l_2 N_t\left(48 + \nicefrac{96 \cdot 27s}{b_x}\right)+ l_1 N_t(24 + 96 \cdot 27s)$ \\
  & \Call{SM-2D}{} & $m + l_2 N_t\left(48 + \nicefrac{96 \cdot 27s}{b_x}\right)+ l_1 N_t(24b_y + 96 \cdot 27s)$ \\
\bottomrule
\end{tabular}
\end{table}

\begin{table}[!h]
\centering
\caption{P2P empirical and analytical time (in milliseconds) on the NVIDIA H200 GPU. $T$ denotes the empirical time. $T_0$, $T_\infty$, $T_{\infty, l_1}$ and $T_f$ denote the analytical time under the 0-L2 cache, $\infty$-L2 cache, $\infty$-L2 cache with prefetching and compute-only model, respectively.  
$l_2 = \nicefrac{1}{4}$, $l_1 = \nicefrac{1}{10}$. $N_t = 4\times10^6$, $s = 256$. }
\label{tab:p2p_time_model}
\setlength{\aboverulesep}{0pt}
\setlength{\belowrulesep}{0pt}
\setlength{\tabcolsep}{3pt}
\begin{tabular}{ccc>{\columncolor[gray]{0.9}}cccc>{\columncolor[gray]{0.9}}cc}
\toprule
& method & $(b_x, b_y)$ &  $T$ & $T_0$ & $T_{\infty}$ &  $T_{\infty, l_1}$ & $T_f$ \\
\midrule
\multirow{4}{*}{P2P}
& \Call{GM-1D}{}& $(64, 1)$& $261$& $854$& $357$& $257$& $190$&\\
& \Call{GM-2D}{}& $(8, 32)$& $336$& $855$& $357$& $257$& $190$&\\
& \Call{SM-1D}{}& $(32, 1)$& $258$& $278$& $262$& $259$& $190$&\\
& \Call{SM-2D}{}& $(64, 2)$& $323$& $267$& $260$& $258$& $190$&\\ 
\bottomrule
\end{tabular}
\end{table}

\subsection{P2G}\label{sec:perf-p2g}
In the P2G kernel, each source interacts with $P^3$ grid points, except for the \Call{P2G-grid}{} method, which interacts with each grid point in colleague cells.
In the \Call{P2G-base}{}, \Call{P2G-source}{}, and \Call{P2G-hybrid}{} methods, $\widetilde{w}_0$ is calculated per source, and the cost can be neglected;
in the \Call{P2G-grid}{} method, $\widetilde{w}_0$ needs to be calculated at each interaction.
The flop count for \Call{P2G-grid}{} is $f_\text{P2G} = 6\nu N_sP^3+ 50\cdot 27 N_s \left(\nicefrac{P}{2}\right)^3$, whereas for the remaining variants
is $f_\text{P2G}=  46 N_s P^3$,

P2G load operations involve source points coordinates, normal vectors, and both single layer and double layer densities ($12N_s$ doubles). P2G stores in memory the vector and tensor densities ($12N_g$ doubles) at the grid points.
Therefore, $m = 96N_s + 96N_g$.
Let $s_f$ be the number of sources in one far-field cell.
In the \Call{P2G-hybrid}{} method, each block stores $b_x \leq s_f$ sources in shared memory, and atomically writes to all grid points in colleague cells.
Thus, we have $m_G^\mathrm{\Call{hybrid}{}} = 96N_s + ( \nicefrac{96\cdot 27 N_s \left ( \nicefrac{P}{2} \right)^3}{b_x} )$. 
The block loads each source's anchors from shared memory for each grid point in colleague cells, and if the grid point is in the $P^3$ cube of the source, it loads the source's normals, two densities, and its approximated window function values.
Therefore, $m_S^\mathrm{\Call{hybrid}{}} = 96 N_s {P}^3 + 12\cdot 27 N_s \left( \nicefrac{P}{2}\right)^3$.
Shared memory is not used for other methods.
For \Call{base}{} and \Call{source}{}, $m_G^\mathrm{\Call{base}{}}  = m_G^\mathrm{\Call{source}{}} = 24N_s + 168 N_s {P}^3$; for \Call{grid}{},
$m_G^\mathrm{\Call{grid}{}} = 96\cdot 27 N_s \left( \nicefrac{P}{2} \right)^3 + 96 N_g$, and $m_S^\mathrm{\Call{base}{}}  = m_S^\mathrm{\Call{source}{}} = m_S^\mathrm{\Call{grid}{}} = 0$. 
\Cref{tab:p2g_mop_all} summarizes the
memory operations.
\Cref{tab:p2g_time_model} shows the empirical and analytical time under different performance models for the P2G kernel.

\begin{table}
\centering
\caption{Analytical MOP values (in bytes) for P2G variants for the $T_{\infty}$ model; $b_x$ is the GPU block size and  $m=96N_s + 96N_g$. The MOP values for the $T_0$ model ($m_0$) are obtained by setting $m = 0$ and $l_2 = 1$, while the MOP values for the $T_{\infty, l_1}$ model ($m_{\infty, l_1}$) are obtained by setting $l_2 = l_1$.}
\label{tab:p2g_mop_all}
\setlength{\aboverulesep}{0pt}
\setlength{\belowrulesep}{0pt}
\setlength{\tabcolsep}{3pt}
\begin{tabular}{lll}
\toprule
\multirow{4}{*}{$m_{\infty}$}
  & \Call{base}{} \& \Call{source}{} &$m + l_2 (24N_s + 168 N_s {P}^3)$ \\
  & \Call{grid}{} & $m + l_2 \left(96\cdot 27 N_s \left( \nicefrac{P}{2} \right)^3 + 96 N_g\right)$\\
  & \Call{hybrid}{} & $m + l_2 \left(96N_s +  \nicefrac{96\cdot 27 N_s \left ( \nicefrac{P}{2} \right)^3}{b_x}\right)$ \\
  &                 & $+\, l_1 \left(96 N_s {P}^3 + 12\cdot 27 N_s \left( \nicefrac{P}{2}\right)^3\right)$ \\
\bottomrule
\end{tabular}
\end{table}

\begin{table}[!h]
\centering
\caption{P2G and G2P empirical and analytical runtime (in milliseconds) on the NVIDIA H200 GPU. $T$ denotes the empirical time. $T_0$, $T_\infty$, $T_{\infty, l_1}$ and $T_f$ denote the analytical time under the 0-L2, $\infty$-L2, $\infty$-L2 with prefetching and compute-only model, respectively.  
$l_2 = \nicefrac{1}{4}$, $l_1 = \nicefrac{1}{10}$. $N_s = 4\times10^6$, $P = 8$, $s = 224$, $\epsilon = 10^{-9}$.}
\label{tab:p2g_time_model}
\setlength{\aboverulesep}{0pt}
\setlength{\belowrulesep}{0pt}
\setlength{\tabcolsep}{3pt}
\begin{tabular}{llr>{\columncolor[gray]{0.9}}cccc>{\columncolor[gray]{0.9}}cc}
\toprule
& method & $(b_x, b_y)$ & $T$ & $T_0$ & $T_{\infty}$ &  $T_{\infty, l_1}$ &$T_f$ \\
\midrule
\multirow{4}{*}{P2G}
& \Call{base}{}&(256, 1)& $477$&$88.9$& $24.9$& $12.0$& $2.81$&\\
& \Call{source}{}&(512, 1)& $354$&$88.9$& $24.9$& $12.0$& $2.81$&\\
& \Call{grid}{}&(128, 1)& $246$&$98.2$& $34.2$& $21.3$& $12.2$&\\
& \Call{hybrid}{}&(64, 1)& $39.5$&$21.7$& $13.4$& $11.6$& $2.81$&\\
\midrule
\multirow{2}{*}{G2P}
& \Call{base}{}&(32, 1) &62.8&12.8& 3.73& 1.88& 0.489&\\
& \Call{target}{}&(32, 1) &28.3&12.8& 3.73& 1.88& 0.489&\\
\bottomrule 
\end{tabular}
\end{table}

\subsection{G2P}\label{sec:perf-g2p}
In both the \Call{G2P-base}{} and \Call{G2P-target}{} methods of the G2P kernel, each thread globally reads one target and $P^3$ grid points, computes the interactions and summations, and then globally writes to one target. Therefore, model parameters are given by
$f = 8N_tP^3$,
$m_0 =  24 N_tP^3$
and 
$m_{\infty}  = 48N_t + 24N_g + l_2 ( 24 N_tP^3)$.

\cref{tab:p2g_time_model} compares predicted runtimes of the P2G and G2P models with empirical results.
The measured times are lower than those predicted by the $T_0$ model.
We speculate that the most likely cause of weak performance is due to uncoalesced accesses of loads and stores. 
The cache provides little benefit due to the irregular memory accesses. 

%% file: sections/results.tex
\section{Numerical Results and Performance Analysis}
\label{sec:results}
The accuracy of the spectral Ewald method is well studied (see \cite{Bagge_Tornberg_2023} for Stokes, \cite{Shamshirgar_2021} for Poisson).
Here, we study the HPC performance of our implementation.

We summarize the key questions of our investigation:
\begin{itemize}[topsep=0pt,itemsep=-1ex,partopsep=1ex,parsep=1ex]
    \item {\bf Performance comparison of the Ewald P2P with simpler P2P kernels:} How does the Ewald Stokes P2P kernel compare to that of free-space Poisson and Stokeslets P2P kernels? How does differing arithmetic intensity affect performance? (\cref{fig:p2p-workload}) 
    \item {\bf Algorithmic variants:} How do the different variants of P2P, P2G, G2P perform on a single device? (\cref{tab:p2p-timing}, \cref{tab:p2g-timings} and \cref{tab:g2p}, respectively)
    \item {\bf Performance portability:} What is the performance portability of our implementation? (\cref{fig:p2p-port} for P2P, \cref{fig:p2g-port} for P2G, \cref{fig:se-port} for the overall Ewald sum)
    \item {\bf Overall Ewald sum on a single device:} What is the overall performance for different particle distributions and at different precisions? How does the sum perform and scale on different devices? (\cref{tab:se-dist}, \cref{tab:se-precision} and \cref{fig:se-port}, respectively)
    \item {\bf Overall Ewald sum on multiple devices using MPI:}  What are the MPI overheads for our solver with multi-device computations? (\cref{tab:mpi-results})
\end{itemize}

{\bf Experimental setup:}
Our experiments were carried out at the Texas Advanced Computing Center (TACC) and the OACISS Center at the University of Oregon. On TACC, we used Lonestar 6 for the A100s and AMD Epyc 7763, and Vista for the H200 and Grace systems; On OACISS we used the Odyssey M1300A system.

The \pylib{} library is implemented in Python, combining CuPy \cite{cupy_learningsys2017} for array operations with PyKokkos \cite{pyKokkos} to interface with the Kokkos kernels described in \cref{sec:acceleration}. The CPU backend runs via the \texttt{Kokkos::OpenMP} execution space using \texttt{gcc} with the \texttt{-O3}, \texttt{-march=native}, and \texttt{-mtune=native} compiler flags, but without explicit vectorization or tuning.
GPU code is compiled via the \texttt{Kokkos::CUDA} execution space for NVIDIA devices and the \texttt{Kokkos::HIP} execution space for AMD devices.
We profile performance using Tau2 \cite{tau} for system-level data and NVIDIA Nsight Compute \cite{nsight-compute} for kernel-level analysis on NVIDIA GPUs.
Unless otherwise noted, source and target $\y_j, \x_i$, as well as $\f(\y_j), \q(\y_j)$ and $\n(\y_j)$, are sampled from a uniform probability distribution in the computational domain with $\Nt=\Ns$
and all code is compiled in double precision.
Flops are counted using the flop-equivalences  in \cref{tab:cycles}.

\subsection{P2P}
\label{p2p:results}

\Cref{tab:p2p-timing} compares execution times and efficiencies for the P2P \Call{GM-1D}{}, \Call{GM-2D}{}, \Call{SM-1D}{}, and \Call{SM-2D}{} methods on a single NVIDIA H200 GPU for different $s$ (number of points per cell) with $N = 10^6$. 
We use an empirically computed optimal block size for each method.
For 1D methods, $b_x$ is selected from $B=\{32, 64, 128, 256, 512\}$; for 2D methods, $(b_x, b_y)$ is selected from $\{(b_x, b_y): b_x b_y\in B,\, b_y\in \{2, 8, 32\}\}$. The source chunk size $b_xb_y$ is constrained by the shared memory size on H200.
For all methods, we achieve $63\%\text{--}78\%$ flop efficiency with $s = 512$. 
The \Call{GM-1D}{} and \Call{SM-1D}{} methods exhibit similar execution time and efficiency, and significantly outperform their 2D counterparts. 
Interestingly, the \Call{GM-2D}{} method performs best when the source dimension is highly parallelized ($b_y = 32$), while the \Call{SM-2D}{} method performs best with near-1D layouts ($b_y = 2$). 
In \Call{SM-2D}{}, a large $b_y$  results in a shared memory allocation of $96 \Nt b_y$, which may reduce the number of blocks launched concurrently due to shared memory size constraints. 

\begin{table}
\centering
\caption{Particle-to-particle (P2P) performance results on a single NVIDIA H200 GPU. Flop efficiency is based on the $T_{f}$ model. $s = 256$. $(b_x, b_y)$ is the optimal.}
\label{tab:p2p-timing}
\setlength{\aboverulesep}{0pt}
\setlength{\belowrulesep}{0pt}
\setlength{\tabcolsep}{3pt}
\begin{tabular}{rr>{\columncolor[gray]{0.9}}c>{\columncolor[gray]{0.9}}c>{\columncolor[gray]{0.9}}ccccc}
\toprule
\multirow{2}{*}{$10^6 \Nt$} & \multirow{2}{*}{method} & \multicolumn{3}{c}{$s = 256$} & \multicolumn{3}{c}{$s = 512$} \\
\cmidrule(lr){3-5} \cmidrule(lr){6-8}
 & & $(b_x,b_y)$ & $\frac{\Nt}{\si{\micro\second}}$ & Flops & $(b_x,b_y)$ & $\frac{\Nt}{\si{\micro\second}}$ & Flops \\
\midrule
1 & GM-1D& $(64, 1)$& 14.8& $71\%$& $(64, 1)$& 8.0& $76\%$&\\
          & GM-2D& $(1, 32)$& 11.8& $56\%$& $(2, 32)$& 6.7& $63\%$&\\
          & SM-1D& $(32, 1)$& 15.2& $72\%$& $(32, 1)$& 8.2& $78\%$&\\
          & SM-2D& $(64, 2)$& 12.4& $59\%$& $(128, 2)$& 6.6& $63\%$&\\
\bottomrule
\end{tabular}
\end{table}

\begin{figure}
    \centering
    \includegraphics[width=0.9\linewidth]{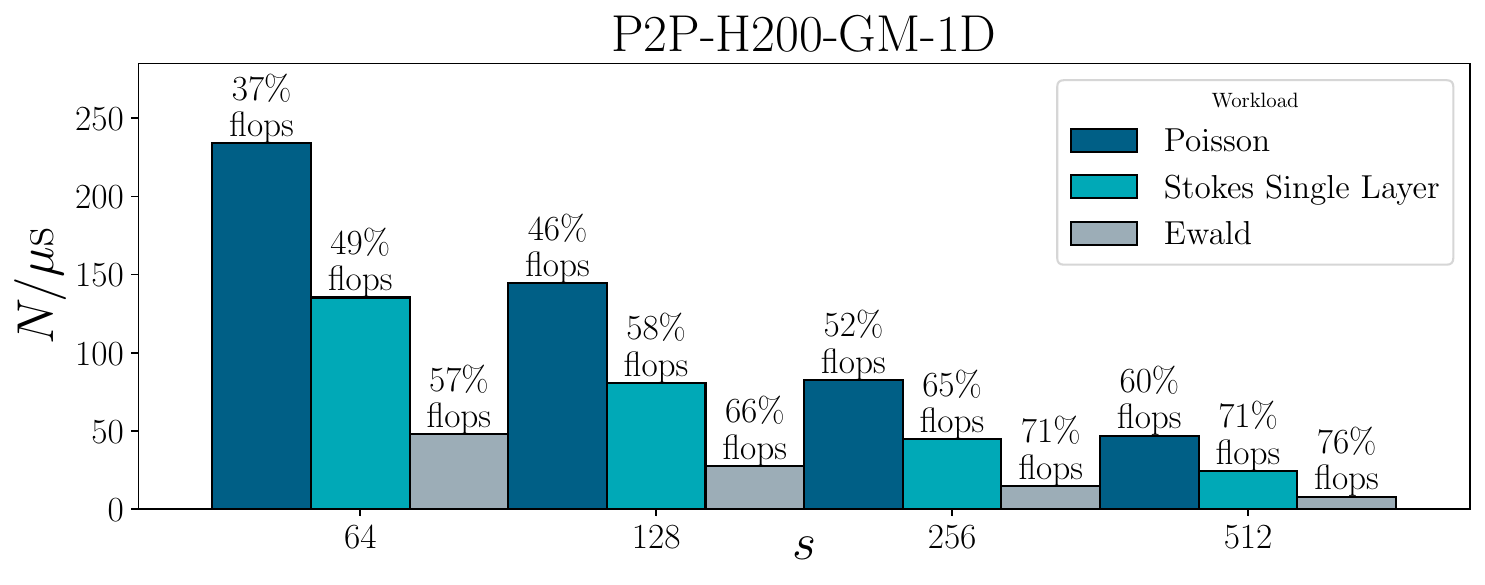}
    \caption{Particle-to-particle performance on a single NVIDIA H200 GPU for different workloads (higher is better), $N = 10^6$.}
    \label{fig:p2p-workload}
\end{figure}

\Cref{fig:p2p-workload} shows the performance of the \Call{P2P-GM-1D}{} kernel on a single NVIDIA H200 GPU under varying workloads, including the Poisson kernel (lightest), the Stokes single-layer kernel, and the full Ewald kernel (heaviest), with $N = 10^6$. For each workload, performance improves with increasing $s$ due to better latency hiding, as each thread processes more interactions.  
The full Ewald kernel achieves the highest efficiency ($76\%$ flops with $s = 512$) as the evaluation of special functions makes it more computationally intensive. Even the lightest workload, the Poisson kernel, achieves $60\%$ flops when $s = 512$, showing that the strong performance of the P2P implementation does not rely on the particular workload and generalizes to other problems.

\begin{figure}
    \centering
    \includegraphics[width=0.9\linewidth]{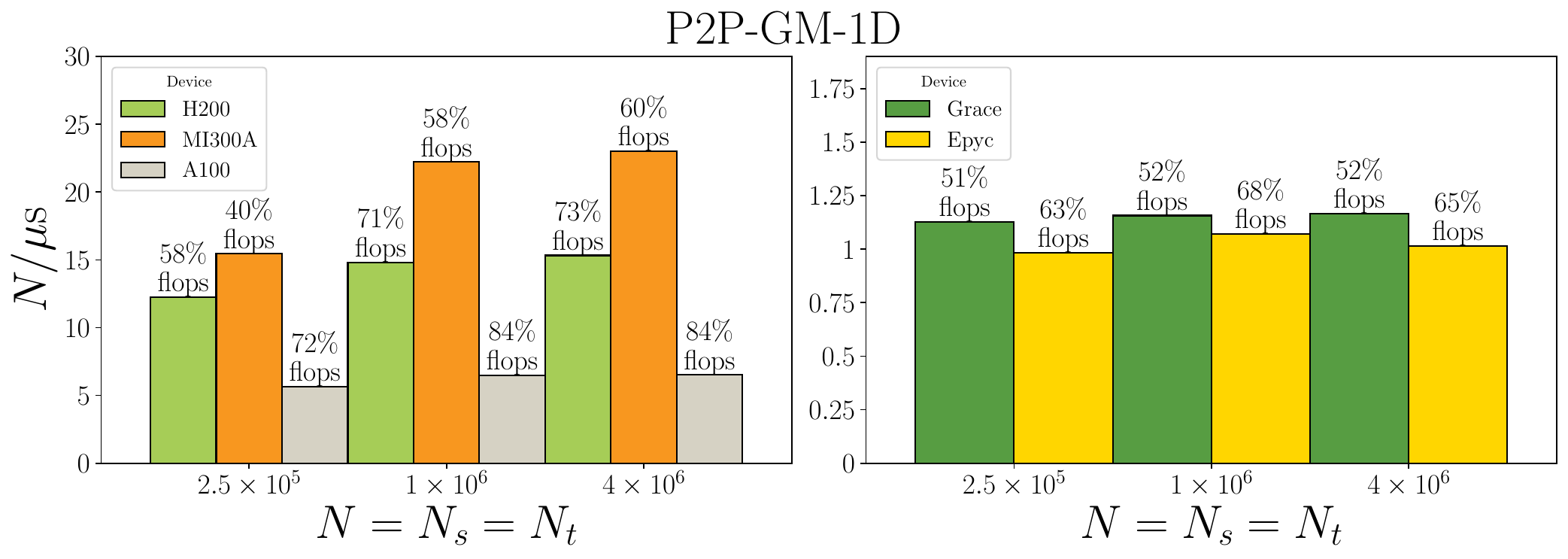}
    \caption{Particle-to-particle performance portability (higher is better), $s = 256$. Note the different $y$-axes for GPU and CPU plots.}
    \label{fig:p2p-port}
\end{figure}

\Cref{fig:p2p-port} presents a performance comparison of the \Call{P2P-GM-1D}{} kernel for the full Ewald sum on various GPUs and CPUs. 
For GPUs,
we achieve flop efficiencies of $84\%$ (A100), $73\%$ (H200), and $60\%$ (MI300A) with $N = 4\times10^6$, $s = 256$. 
NVIDIA devices achieve higher efficiency, but the AMD MI300A attains the lowest execution time due to its highest peak flop performance. On CPUs, the NVIDIA Grace achieves over $50\%$ efficiency, while AMD EPYC exceeds $60\%$. Efficacy is achieved with the ported Kokkos kernels \emph{as is}, without hardware-specific tuning.

\subsection{P2G}
\label{p2g-results}

\Cref{tab:p2g-timings} compares different P2G methods with varying numbers of source particles and target accuracies $E$.
Each row reports results using an empirically computed optimal 1D thread configuration, $b_x \in \{32, 64, 128, 256, 512\}$, on a single NVIDIA H200 GPU.
While \Call{P2G-source}{}, \Call{P2G-grid}{}, and \Call{P2G-hybrid}{} maintain consistent efficiency for fixed accuracy and increasing $\Ns$, the baseline \Call{P2G-base}{} method \emph{loses} efficiency as $\Ns$ increases, resulting in larger speedups for \Call{P2G-hybrid}{}.
This trend probably results from increasing atomic conflicts in \Call{P2G-base}{} with larger problem sizes.
Interestingly, \Call{P2G-source}{} underperforms relative to \Call{P2G-base}{} at low accuracies, which we attribute to increased atomic conflicts due to increased spatial locality within thread blocks.

\begin{table}
\centering
\caption{Particle-to-grid (P2G) performance results on a single NVIDIA H200 GPU. Efficiency is based on the $T_0$ model. Near-field cell size $s$ is $160$ for $P=4, E=10^{-7}$ and $224$ for $P=8, E=10^{-9}$. Each row is sub-selected as the fastest run with $b_x \in \{ 32, 64, 128, 256, 512 \}$. Efficiency (Eff) is reported in mops for the BASE, SOURCE, and HYBRID methods, and flops for the GRID method.}
\label{tab:p2g-timings}
\setlength{\aboverulesep}{0pt}
\setlength{\belowrulesep}{0pt}
\setlength{\tabcolsep}{3pt}
\begin{tabular}{rl>{\columncolor[gray]{0.9}}c>{\columncolor[gray]{0.9}}c>{\columncolor[gray]{0.9}}rccr>{\columncolor[gray]{0.9}}c>{\columncolor[gray]{0.9}}c>{\columncolor[gray]{0.9}}r}\toprule
        \multirow{2}{*}{$10^6 \Ns$} 
        & \multirow{2}{*}{method} 
        &\multicolumn{3}{c}{$P=4,E=10^{-7}$}
        &\multicolumn{3}{c}{$P=8, E=10^{-9}$}
        &\multicolumn{3}{c}{$P=14, E=10^{-14}$}\\
	\cmidrule(r){3-5}\cmidrule(r){6-8}\cmidrule{9-11}
         &&$\frac{\Ns}{\si{\micro\second}}$&Speedup&Eff            
	 &$\frac{\Ns}{\si{\micro\second}}$&Speedup&Eff 
	 &$\frac{\Ns}{\si{\micro\second}}$&Speedup&Eff \\
\midrule
	     $0.25$& \Call{base}{}  &$101.65$ & ---          &$6\%$  &$12.52$  & ---           &$5\%$  & $1.91$  & ---           & $3\%$ \\
		   &\Call{source}{} &$48.2$   & $0.47\times$ &$3\%$  &$11.83$  & $0.95\times$  &$5\%$  & $2.5 $  & $1.31\times$  & $3\%$ \\
		   &\Call{grid}{}   &$83.41$  & $0.82\times$ &$10\%$ &$16.42$  & $1.31\times$  &$26\%$ & $3.14$  & $1.64\times$  & $49\%$\\
		   &\Call{hybrid}{} &$147.22$ & $1.45\times$ &$9\%$  &$90.73$  & $7.25\times$  &$37\%$ & $21.93$ & $11.48\times$ & $31\%$\\
\midrule
	     $1$   & \Call{base}{}  &$103.02$ & ---          &$6\%$  &$9.97$   & ---           &$4\%$  & $1.61$  & ---           & $2\%$\\
		   &\Call{source}{} &$62.36$  & $0.61\times$ &$4\%$  &$11.94$  & $1.20\times$  &$5\%$  & $2.4$   & $1.49\times$  & $3\%$\\
		   &\Call{grid}{}   &$93.58$  & $0.91\times$ &$11\%$ &$16.87$  & $1.69\times$  &$27\%$ & $3.1$   & $1.92\times$  & $48\%$\\
		   &\Call{hybrid}{} &$177.26$ & $1.72\times$ &$11\%$ &$100.65$ & $10.09\times$ &$42\%$ & $22.1$  & $13.72\times$ & $27\%$\\
\midrule
	     $4$   & \Call{base}{}  &$80.57$  & ---          &$5\%$  &$8.38$   & ---           &$3\%$  & $1.49$  & ---           & $2\%$\\
		   &\Call{source}{} &$54.17$  & $0.67\times$ &$3\%$  &$11.3$   & $1.35\times$  &$4\%$  & $2.52$  & $1.69\times$  & $3\%$\\
		   &\Call{grid}{}   &$82.52$  & $1.02\times$ &$10\%$ &$16.28$  & $1.94\times$  &$26\%$ & $3.39$  & $2.28\times$  & $52\%$\\
		   &\Call{hybrid}{} &$178.23$ & $2.21\times$ &$11\%$ &$101.28$ & $12.09\times$ &$39\%$ & $24.05$ & $16.18\times$ & $25\%$\\
\bottomrule

        \end{tabular}
    \end{table}

At high accuracy tolerance, the advantage of \Call{P2G-hybrid}{} becomes more pronounced.
For example, with $\Ns = 4 \times 10^6$, spectral Ewald tolerance $10^{-14}$ and optimal cell size $s = 960$, \Call{P2G-hybrid}{} achieves a \textbf{16.18$\times$} speedup 
over \Call{P2G-base}{}.
In this regime, the large near-field $s$ leads to a smaller Fourier grid, which again increases atomic contention in \Call{P2G-base}{} while conversely enhancing spatial locality in \Call{P2G-hybrid}{}.

\begin{figure}
    \centering
    \includegraphics[width=0.9\linewidth]{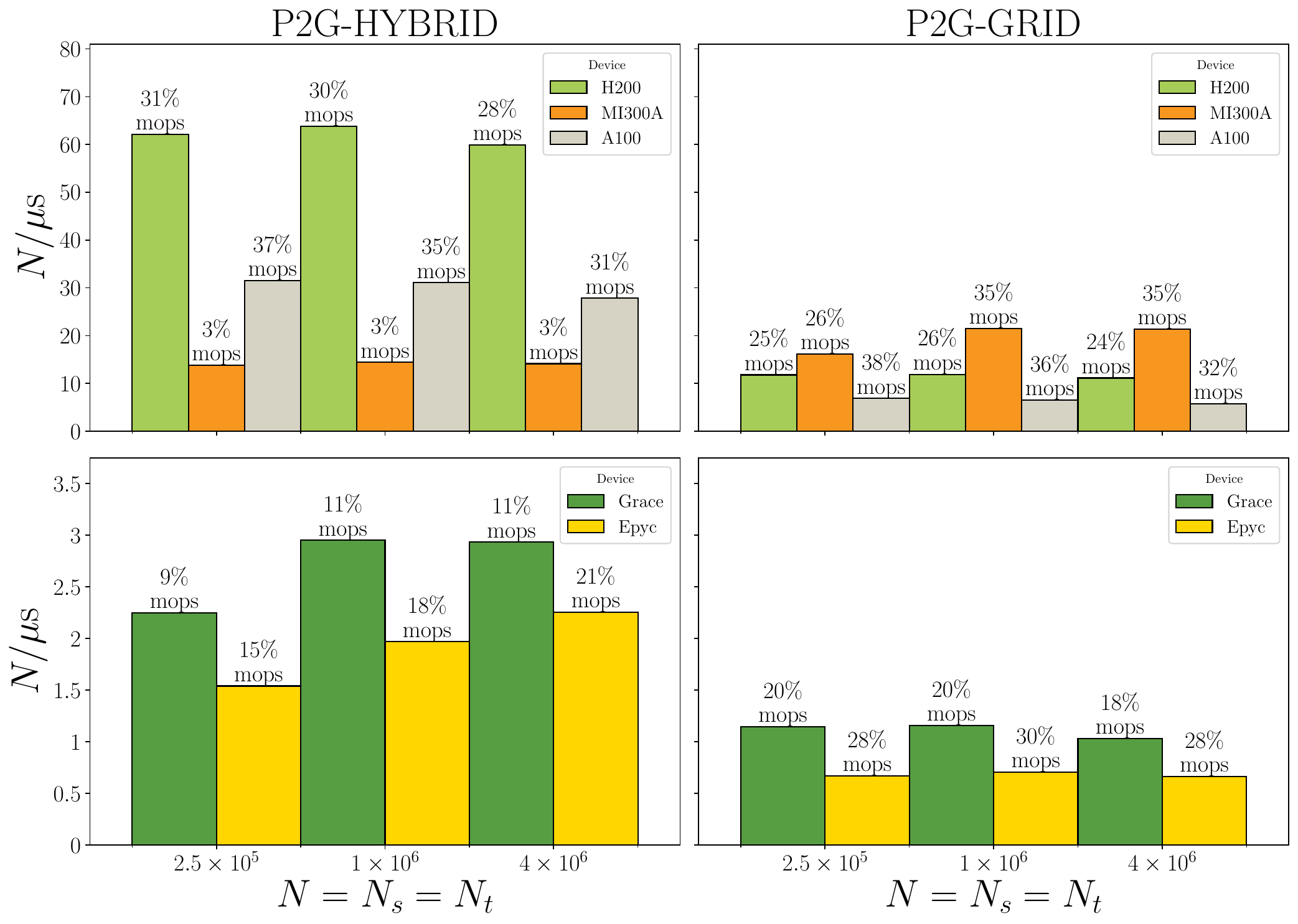}
    \caption{Particle-to-grid performance portability (higher is better). $E = 10^{-9}$, $s = 224$. Note the different $y$-axes for GPU and CPU plots.}
    \label{fig:p2g-port}
\end{figure}

\Cref{fig:p2g-port} explores the performance portability of \Call{P2G-hybrid}{} and \Call{P2G-grid}{}. The fastest overall performance is achieved by \Call{P2G-hybrid}{}, which processes approximately $6 \times 10^7$ particles per second with $30\%$ efficiency. Surprisingly, the AMD MI300A underperforms significantly for the \Call{P2G-hybrid}{} method, processing nearly \emph{6$\times$ fewer} particles per second than NVIDIA GPUs.
This drop appears to be related to shared memory limitations in AMD devices, consistent with the trends observed for \Call{P2P-SM-1D}{}. \Call{P2G-hybrid}{} has a shared memory workload of $(3P + 9)\frac{s_f}{b_x}$.  On the other hand, \Call{P2G-grid}{} exhibits consistent performance and efficiency on both NVIDIA and AMD GPUs, and in fact outperforms \Call{P2G-hybrid}{} on the MI300A. On CPUs, \Call{P2G-hybrid}{} remains the optimal method for P2G computations, although its performance advantage and achieved efficiencies are less pronounced compared to the GPUs. 

\subsection{G2P}
\label{g2p-resultss}

\Cref{tab:g2p} compares \Call{G2P-base}{} and \Call{G2P-target}{} for varying numbers of target particles $\Nt$ and spectral Ewald tolerances (controlled by $P = 4$ and $P = 8$) on a single NVIDIA H200 GPU.
In general, low-accuracy computations ($P = 4$) are less efficient.
For smaller problem sizes, \Call{G2P-base}{} outperforms \Call{G2P-target}{} due to lower GPU overhead and simpler memory access patterns.

\begin{table}
\centering
        \caption{Grid-to-particle (G2P) performance results on a single NVIDIA H200 GPU. Mop efficiency is calculated using the $T_0$ model. Near-field cell size $s$ is $160$ for $P=4, E=10^{-7}$ and $224$ for $P=8, E=10^{-9}$. Each row is sub-selected as the fastest run with $b_x\in \{ 32, 64, 128, 256, 512 \}$.}\label{tab:g2p}
\setlength{\aboverulesep}{0pt}
\setlength{\belowrulesep}{0pt}
\setlength{\tabcolsep}{3pt}
\begin{tabular}{rl>{\columncolor[gray]{0.9}}r>{\columncolor[gray]{0.9}}c>{\columncolor[gray]{0.9}}rrcr}\toprule
        \multirow{2}{*}{$10^6\Nt$} 
        & \multirow{2}{*}{method} 
        &\multicolumn{3}{c}{$P=4, E=10^{-7}$}
        &\multicolumn{3}{c}{$P=8, E=10^{-9}$}\\
        \cmidrule(r){3-5}\cmidrule(r){6-8}
         &&$\frac{\Nt}{\si{\micro\second}}$&Speedup&Mops            &$\frac{\Nt}{\si{\micro\second}}$&Speedup&Mops \\
        \cmidrule(r){3-5}\cmidrule(r){6-8}

\midrule
$0.25$& \Call{base}{}&$447.73$& ---&$21\%$&$85.56$& ---&$31\%$\\
   & \Call{target}{}&$430.85$& $0.96\times$&$20\%$&$119.41$& $1.40\times$&$44\%$\\
\midrule
$1$& \Call{base}{}&$661.48$& ---&$28\%$&$66.53$& ---&$24\%$\\
   & \Call{target}{}&$668.02$& $1.09\times$&$31\%$&$136.97$& $2.06\times$&$50\%$\\
\midrule
$4$& \Call{base}{}&$523.13$& ---&$24\%$&$63.68$& ---&$23\%$\\
   & \Call{target}{}&$760.22$& $1.45\times$&$35\%$&$141.38$& $2.22\times$&$52\%$\\
\bottomrule

        \end{tabular}
\end{table}

At higher accuracy ($P = 8$), \Call{G2P-target}{} becomes more efficient, as expected.
This improvement is primarily due to better memory coalescing over the target particles, enabled by the introduction of the cell list $\CF$.

\Cref{fig:se-port} shows the performance portability of \Call{G2P-target}{} in the context of a full Ewald sum with $s = 224$ and target error $E = 10^{-9}$. G2P efficiency is modeled in MOP/s using the 0-cache model. As has become a pattern with our memory-bound Kokkos kernels, we observe similar efficiencies on NVIDIA GPUs, while the AMD MI300A exhibits significantly lower efficiency. For the same run on CPUs, i.e., $s=224$, $N_s = 2.5\times 10^5, 1\times 10^6, 4\times 10^6$, we observe efficiencies of 50\%, 60\%, and 58\% on NVIDIA Grace and 65\%, 86\%, and 93\% on AMD Epyc.

\subsection{Overall Ewald Sum}
\label{sec:se-results}

\cref{tab:se-dist} compares the run-times for each stage in the Ewald sum for different particle distributions with target error $E=10^{-9}$, near-field cell size $224$, and $N=4\times 10^6$ source and target points. 
The P2P and P2G stages are most sensitive to non-uniform particle distribution as both algorithms are reliant on \emph{non-adaptive} cell lists that grow in size are particles become clustered.
Detailed analysis of Ewald sum runtimes for particle distributions with varying levels of non-uniformity are provided in supplement \ref{sec:nonuniform}. 

\begin{table}
\label{tab:se-dist}
\centering
\caption{Ewald sum performance in milliseconds on a single NVIDIA H200 GPU for different particle distributions; $N=4\times 10^6$, $s=224$, $P=8$ and $E=10^{-9}$. \textbf{X/Y} indicates the source and target distribution, respectively, where \textbf{U}(0,1) is a uniform-distribution, \textbf{N}(0,0.3) a normal-distribution shifted and truncated to the unit cube, and \textbf{S} is uniform on the surface of the unit sphere. We use the P2P-GM-1D, P2G-HYBRID, and G2P-TARGET methods.}
\setlength{\aboverulesep}{0pt}
\setlength{\belowrulesep}{0pt}
\setlength{\tabcolsep}{3pt}
\begin{tabular}{c|ccccc}\toprule
& \U/\U & \N/\N & \N/\U & \S/\S & \S/\U\\
\midrule
P2P &239.0 & 351.4 & 249.4 & 901.8 & 268.9\\
P2G &47.6 & 77.4 & 76.3 & 59.7 & 59.2\\
FFT &48.7 & 48.7 & 48.7 & 48.8 & 48.8\\
CNV &20.0 & 20.1 & 20.2 & 20.6 & 20.4\\
IFFT &10.7 & 10.7 & 10.7 & 10.8 & 10.8\\
G2P &36.7 & 37.1 & 36.6 & 21.3 & 36.8\\
\midrule
Total &402.8 & 545.4 & 441.9 & 1063.0 & 444.8\\
\bottomrule
    \end{tabular}
\end{table}

\cref{tab:se-precision} compares the runtimes of the Ewald sum in single and double precision arithmetic with target error $E=10^{-7}$, near-field cell size $160$, and $N=4\times 10^6$ source and target points. 
All Ewald stages \emph{besides} P2G and IFFT show about a $2\times$ speedup when performed in single precision arithmetic. 
The P2G algorithm involves the same costly integer arithmetic and block synchronization in both single and double floating point precision, yielding only a $1.3\times$ speedup. 
We hypothesize that the IFFT launch is dominated by overhead for such small a size in both precisions. 
Overall, we see a $1.9\times$ speedup between the single precision and double precision Ewald sums. 

\begin{table}
\label{tab:se-precision}
\centering
\caption{Spectral Ewald performance in milliseconds with $N=4\times10^6$, results on a single NVIDIA H200 GPU for different dtypes. Near-field cell size $s$ is $160$ and the target error $E=10^{-7}$. We use the P2P-GM-1D, P2G-HYBRID, and G2P-TARGET methods.}
\setlength{\aboverulesep}{0pt}
\setlength{\belowrulesep}{0pt}
\setlength{\tabcolsep}{3pt}
\begin{tabular}{c|ccccc}
\toprule
& single-precision & double-precision\\
\midrule
P2P   &92.67  & 191.34\\
P2G   &28.24  & 36.69\\
FFT   &4.8    & 8.17\\
CNV   &6.72   & 11.72\\
IFFT  &1.27   & 1.87\\
G2P   &7.95   & 17.51\\
Total &141.65 & 267.3\\
\bottomrule

    \end{tabular}
\end{table}

\begin{figure}
    \centering
    \includegraphics[width=\linewidth]{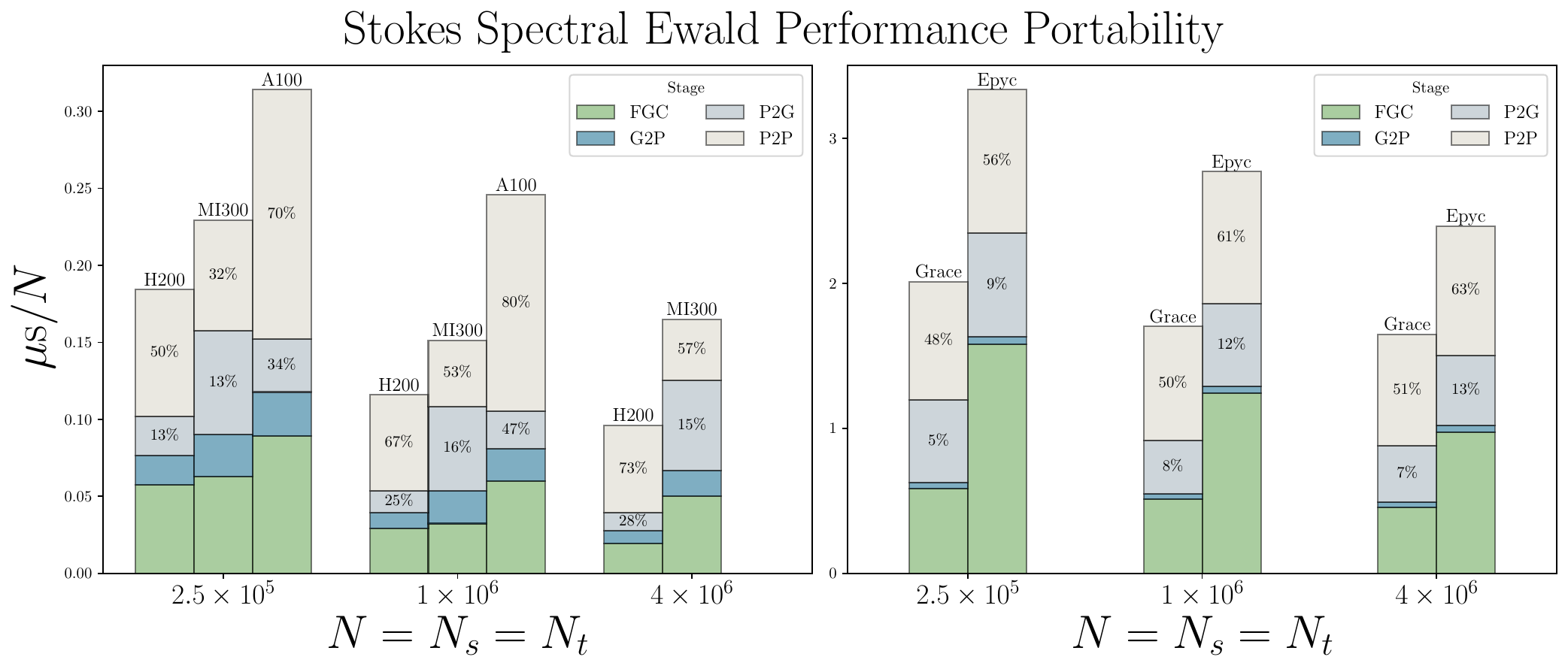}
    \caption{Spectral Ewald run on different machines (lower is better). $E = 10^{-9}$, $s = 224$. P2P efficiencies are reported with the L2 prefetching model (i.e., flop dominant) and the P2G and G2P efficiencies are computed with the 0-cache memory model (i.e., mop dominant). Note the different $y$-axes for GPU and CPU plots.}
    \label{fig:se-port}
\end{figure}

\Cref{fig:se-port} presents a full Ewald sum with fixed target error $E = 10^{-9}$ and near-field cell size $s = 224$. The largest run, with $\Ns = 4 \times 10^6$ source particles, requires an upsampled Fourier grid too large to fit in memory on the A100 GPU. On the MI300A, the absence of a highly efficient P2G kernel causes the total runtime to lag behind that of the H200, despite the MI300A’s higher theoretical throughput. Both CPUs have similar runtimes and efficiencies for the P2P, P2G, and G2P steps. 

The Fourier Grid Convolution step --- dominated by memory-bound FFT and IFFT operations --- is most efficiently executed on the H200 for GPUs, reflecting their superior memory bandwidth and optimized FFT performance. 
Grace dominates the FCG step on the CPU, which we attribute to optimized library implementations for the FFT and IFFT steps.

Overall, as the number of particles increases, we observe clear $\mathcal{O}(N \log N)$ scaling in the overall runtime, consistent with the expected complexity of the Ewald sum algorithm.

\subsection{Multi-GPU Performance}
\label{sec:multi}

We consider a standard multi-GPU strategy based on message passing (MPI).
Since a number of the kernels (P2P, P2G, CNV, G2P) are spatially local, we divide the domain $\B$ into equisized slabs along the periodic spatial dimension, delegating each slab to one GPU.
The slab partition is in fact compatible to problems of any periodicity; the only periodicity specific detail being the edge-slab communication pattern.
Communication between GPUs can then be reduced to: 
\begin{enumerate*}[label=(\textbf{\arabic*})]
    \item initial radix sorting of sources and targets so that each particle belongs to the correct GPU ({\small MPI-SORT});
    \item exchanging ghost source points with neighbors; the points within distance $\rc$ to the slab interface are exchanged before P2P ({\small MPI-GHOST-SOURCE});
    \item similar exchange of ghost grid values between IFFT and G2P steps ({\small MPI-GHOST-GRID});
    \item FFT and IFFT steps, which we implement using NVIDIA's cuFFTMp library.
\end{enumerate*}
Since Kokkos \emph{does not} provide a performance portable distributed FFT, we rely on vendor libraries; the aforementioned cuFFTMp for NVIDIA GPUs, rocFFT for AMD GPUs, and FFTW MPI for CPUs. 

A weak scaling test with $N = \num{4e6} \NGPU$ source and target particles on $\NGPU$ H200 GPUs is shown in \cref{tab:mpi-results}.
The GPUs are connected via NVIDIA InfiniBand with \SI{400}{Gbit/s} bandwidth. Inter-GPU communication takes place in the FFT, IFFT and MPI steps.
We observe bounded $\mathcal{O}(1)$ communication costs as we increase $\NGPU$ in all steps except {\small MPI-SORT}, which is all-to-all and grows approximately linearly.
We choose to implement {\small MPI-SORT} as an all-to-all operation to show worst-case runtime when we have no \emph{a-priori} information about what particles are on what GPUs; 
in many settings, e.g., timestepping regimes where GPUs are assigned spatially close particles, 
{\small MPI-SORT} can be implemented as a local neighbor communication, reducing communication complexity from $\mathcal{O}(\NGPU)$ to $\mathcal{O}(1)$.
It should be noted that the method parameters (e.g., $\rc$, $\xi$, $h$) were not tuned for each $\NGPU$ case, but rather fixed to their single-GPU values.

\begin{table}
    \centering
    \caption{Runtime in \emph{milliseconds} for multi-GPU weak scaling test with $\NGPU$ NVIDIA H200 GPUs, MPI communication, and NVIDIA's cuFFTMp library. Parameters are as in \cref{fig:se-port} with $N = \num{4e6} \NGPU$. Note that the FFT stage comprises 12 real-to-complex distributed FFTs while the IFFT stage comprises only 3 complex-to-real distributed inverse FFTs.}
    \label{tab:mpi-results}
    \begin{tabular}{c rrrrrrr}
    \toprule
    Step & {$\NGPU=1$} & {2} & {4} & {8} & {16} & {32} & 64 \\
    \midrule
P2P &227.4 &232.2 &232.2 &234.2 &235.1 &231.5 &229.9\\
P2G &45.9 &53.6 &51.3 &51.6 &49.9 &50.5 &45.0\\
FFT  & 48.8 & 190.4 & 245.9 & 278.3 & 299.0 & 307.0 & 330.2\\
CNV &18.1 &25.4 &24.8 &24.6 &24.2 &24.3 &23.4\\
IFFT &  10.7 & 47.4 & 61.5 & 69.6 & 74.6 & 76.6 & 83.2\\
G2P &33.3 &35.2 &34.0 &32.7 &32.0 &31.4 &30.5\\
\midrule
\footnotesize{MPI-SORT} &--- &16.5 &23.7 &37.1 &61.3 &110.4 &222.0\\
\footnotesize{MPI-GHOST-SOURCE} &--- &7.6 &8.1 &8.0 &8.0 &8.0 &7.2\\
\footnotesize{MPI-GHOST-GRID} &--- &6.9 &6.2 &6.1 &5.7 &5.5 &2.7\\
    \bottomrule
    \end{tabular}
\end{table}

The FFT runtime begins to dominate as the number of GPUs grows. On the GH200, inter-GPU communication is the bottleneck, limited by the RDMA GPUDirect off-chip over NDR InfiniBand (400~Gbit/s).
We estimate the communication cost for the FFT stage, which comprises 12 real-to-complex distributed FFTs, as 
$ T_\mathrm{comm} \approx 12 \times \left(1-\nicefrac{1}{N_{\text{GPU}}}\right)
\times \nicefrac{N_g}{2 N_{\text{GPU}}} \times \nicefrac{\SI{16}{B}}{b_\text{comm}}$, where  $b_\text{comm}$ is the bandwidth.
The infiniband setup on TACC's Vista computer has $b_\text{comm} = \SI{50}{GB/s}$ and our test problem parameters yield $N_g = 256 N_\text{GPU} \times 692 \times 692$.
As such, the FFTs achieve $70\%$ of this peak communication model, which is in line with NVIDIA's results \cite{Cambier2022blog}.

Note that there are multiple strategies to improve FFT performance, such as tuning the near-field cell size $s$ to shift some cost into the P2P computations; however, we kept this parameter fixed to highlight communication overheads.
It should be noted that the P2P computations are independent of the FFTs and can be overlapped with them, but we report each step independently to report isolated costs per step.

%% file: sections/conclusion.tex
\section{Conclusions} \label{sec:conclusion}
We introduced algorithms designed for the spectral Stokes Ewald sum, accompanied by a comprehensive analysis and evaluation of their performance. In particular, we compared different algorithmic variants as well as the performance of the method on different hardware and for different particle distributions and precisions.
The MPI results served primarily as a preliminary glimpse, highlighting potential overheads associated with off-node operations.
In future research, we will study the optimization of CPU kernels and thorough performance characterization and optimization of the MPI implementation, particularly in multi-GPU NUMA architectures. We also note that novel Ewald decompositions based on prolate spheroidal wave functions (PSWFs) have the potential to further reduce the computational cost \cite{stokes_pswf} compared to the Hasimoto decomposition that was used here. We plan to introduce such a PSWF decomposition into our algorithm in the future.

%% file: sections/appendix.tex
% Appendix
\appendix

\section{The ParkiPy Package}
\label{sec:parkipy}
The algorithms presented in this paper are implemented in the \emph{(Par)ticle (k)ernel (i)nteractions for (Py)thon} package, \pylib{}.
The package offers easy-to-use and architecture-agnostic APIs to evaluate the Stokes and Poisson kernels in arbitrary periodicity using Ewald summation.
\pylib{}'s Ewald summation submodule, \verb|parkipy.ewald|, currently supports the Ewald summation for the Stokes single layer potential, the combined Stokes single and double layer potential (\cref{fig:code}) and the Poisson kernel. There is also a distributed submodule, \verb|parkipy.distributed.ewald|, which computes the singly periodic combined Stokes potential for multi-node multi-gpu settings in a CUDA execution space. The main namespace has a class \verb|parkipy.CellList| for efficient cell list construction.
\begin{figure}
	\begin{python}
import parkipy
ex = parkipy.utils.get_execution_space("CUDA") # Kokkos backend
am = parkipy.utils.get_array_module(ex) # `cupy' or `numpy', cond. on ex

# set spectral Ewald box and tolerance
box = [1.0, 1.0, 1.0]
tol = 1e-4

# generate sources, targets, densities, and normals
nt = 4000000
ns = 4000000
trg = am.random.rand(3, nt) * am.array(box).reshape(3, 1)
src = am.random.rand(3, ns) * am.array(box).reshape(3, 1)
dens_sl = am.random.randn(3, ns)
dens_dl = am.random.randn(3, ns)
norms = am.random.randn(3, ns)
dens = am.vstack((dens_sl, dens_dl))  # stack densities for ewald call

# get the potential for the Stokes single and double layer potential
options = parkipy.ewald.EwaldOptions(
	periodicity=1, box=box, tolerance=tol, 
	cell_size=224, execution_space=ex
)
pot = parkipy.ewald.stokes_comb(trg, src, dens, norms, options)
\end{python}
\caption{{Code using the \pylib{} library to compute the combined Stokes single and double layer potential. This example runs on NVIDIA GPUs and uses the cupy array module as a backend.}}
\label{fig:code}
\end{figure}

\section{Non-uniform Particle Distributions}
\label{sec:nonuniform}

\begin{figure}
    \centering
    \includegraphics[width=\linewidth]{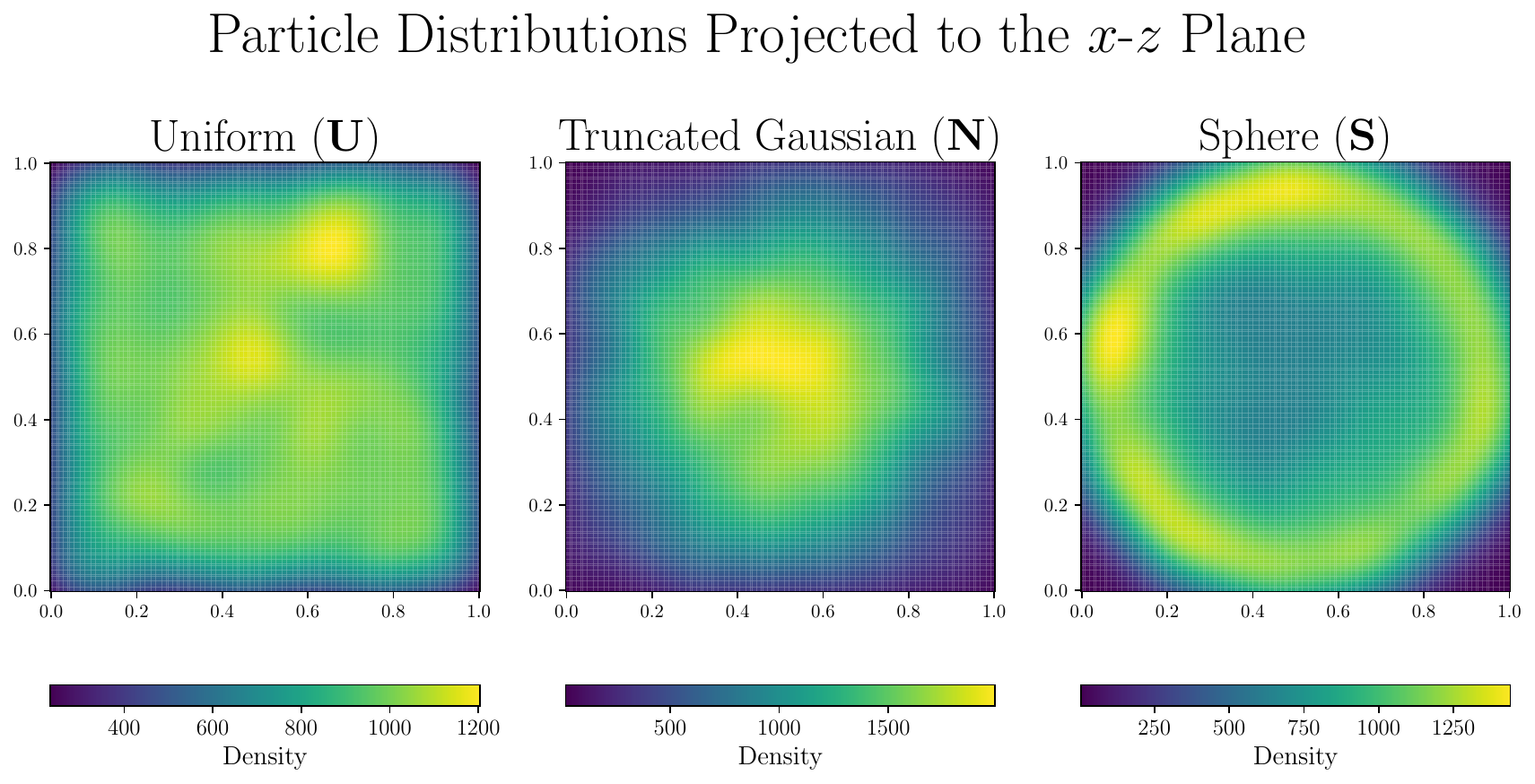}
    \caption{$N=4\times 10^6$ particles projected onto the $x$-$z$ plane. Densities estimated with a Gaussian KDE. \emph{Left}: particles distributed uniformly in the unit cube. \emph{Middle}: particles distributed normally with variance $0.3$. \emph{Right}: particles distributed uniformly on the surface of the unit sphere.}
    \label{fig:dist}
\end{figure}

The performance of the Ewald sum becomes suboptimal for highly non-uniform distributions of the sources and targets. We compare cases where particles are uniformly distributed in the unit cube (\U), normally distributed with variance $0.3$ shifted and truncated to the unit cube (\N), and uniformly distributed on the surface of the unit sphere (\S) (see \cref{fig:dist}). 
When both sources and targets are distributed as \N, our P2P algorithm achieves $68\%$ of the performance relative to a source and target distribution as \U, while our P2G algorithm achieves $61\%$ of the performance.  
When both sources and targets are distributed as \S, our P2P algorithm achieves $27\%$ performance and our P2G algorithm $80\%$ performance, again relative to the performance with sources and targets distributed as \U. 
The P2P algorithm performs close to the ideal distribution of \U when the sources are distributed as \N or \S \emph{but} the targets are distributed as \U, as P2P is parallel over target cell lists and assumes a consistent number of points per target cell.
The G2P algorithm performs best when the targets are distributed as \S, about a $1.7\times$ speedup compared to a distribution of \U, as more targets read the same grid points, increasing cache locality. 
In our worst case, where both sources and targets are distributed as \S, the total Ewald sum achieves $38\%$ of the performance of the best case where sources and targets are distributed as \U. 

